\documentclass[11pt]{amsart}
\usepackage[margin=1in]{geometry}

\usepackage[latin1]{inputenc}
\usepackage[T1]{fontenc}

\usepackage[english]{babel}
\usepackage{lmodern}
\usepackage{amsmath,amssymb,amscd}
\usepackage{fancyhdr}
\usepackage[T1]{fontenc}
\usepackage{amsfonts}
\usepackage{amsthm}
\usepackage{amstext}
\usepackage{a4wide}
\usepackage{dsfont}
\usepackage{enumerate}
\usepackage{exscale}
\usepackage{fancyhdr}
\usepackage{hyperref}
\usepackage{mathtools}
\usepackage{mathrsfs}
\usepackage{paralist}
\usepackage{relsize}
\usepackage{tikz}
\usepackage{upgreek}
\usepackage{url}
\usepackage{verbatim}
\usepackage{graphicx,caption}

\usepackage{appendix}
\usepackage{latexsym}

\newenvironment{paragr}[1][]{\refstepcounter{subsubsection} \noindent \textbf{\thesubsubsection . \ #1}}{\medskip}

\newenvironment{theo}{ \medskip\refstepcounter{thm}  \noindent\textbf{Theorem \thethm}. ---\em}{\em \medskip}
\newenvironment{prop}{\medskip\refstepcounter{thm}   \noindent\textbf{Proposition \thethm}. ---\em}{\em\medskip}
\newenvironment{cor}{\medskip\refstepcounter{thm}  \noindent\textbf{Corollary \thethm}. ---\em}{\em\medskip}

\newenvironment{lem}{\medskip\refstepcounter{thm}   \noindent\textbf{Lemma \thethm}. ---\em}{\em\medskip}
\newenvironment{clm}{\medskip\refstepcounter{thm}   \noindent\textbf{Claim \thethm}. ---\em}{\em\medskip}

\newenvironment{defi}{\medskip\refstepcounter{thm}  \noindent\textbf{Definition \thethm}. ---}{\medskip}

\newenvironment{rem}{\medskip\refstepcounter{thm}  \noindent\textbf{Remark \thethm}. ---}{\medskip}

\makeatletter
\@addtoreset{equation}{section}         
\makeatother

\def\ps#1#2%
  {\mathop{}%
   \mathopen{\vphantom{#2}}^{#1}%
   \kern-\scriptspace%
   #2}

\newcommand{\beqyn}{\begin{eqnarray}}
\newcommand{\enqyn}{\end{eqnarray}}
\newcommand{\blem}{\begin{lem}}
\newcommand{\elem}{\end{lem}}
\newcommand{\brop}{\begin{prop}}
\newcommand{\erop}{\end{prop}}
\newcommand{\bcor}{\begin{cor}}
\newcommand{\ecor}{\end{cor}}
\newcommand{\brem}{\begin{rem}}
\newcommand{\erem}{\end{rem}}
\newcommand{\bclm}{\begin{clm}}
\newcommand{\eclm}{\end{clm}}
\newcommand{\benum}{\begin{enumerate}}
\newcommand{\enum}{\end{enumerate}}
\newcommand{\btheo}{\begin{theo}}
\newcommand{\etheo}{\end{theo}}
\newcommand{\bdem}{\begin{proof}}
\newcommand{\edem}{\end{proof}}
\newcommand{\bdefi}{\begin{defi}}
\newcommand{\edefi}{\end{defi}}
\newcommand{\bpar}{\begin{paragr}}
\newcommand{\epar}{\end{paragr}}

\newcommand{\scrF}{\mathscr{F}}

\newcommand{\Hom}{\mathrm{Hom}}

\newcommand{\Lie}{\mathrm{Lie}}

\newcommand{\Ad}{\mathrm{Ad}}

\newcommand{\End}{\mathrm{End}}

\newcommand{\Res}{\mathrm{Res}}

\newcommand{\Rel}{\mathrm{Re}}

\newcommand{\upla}{\underline{\rho}}
\newcommand{\al}{\alpha}

\newcommand{\la}{\lambda}
\newcommand{\La}{\Lambda}

\newcommand{\ssum}[1]{\sum_{\begin{subarray}{c} #1 \end{subarray}}}

\newcommand{\bilif}{\langle \cdot,\cdot \rangle}

\newcommand{\dsl}{\displaystyle \left(}
\newcommand{\rb}{\right)}
\newcommand{\dsla}{\displaystyle \left|}
\newcommand{\rba}{\right |}
\newcommand{\bsl}{\backslash}

\newcommand{\hrar}{\hookrightarrow}

\newcommand{\smin}{\smallsetminus}
\newcommand{\sbs}{\subset}
\newcommand{\sbn}{\subsetneq}
\newcommand{\nsbs}{\not \subset}
\newcommand{\sps}{\supset}

\newcommand{\ago}{\mathfrak{a}}

\newcommand{\dl}{\mathfrak{d}}

\newcommand{\Ugo}{\mathfrak{U}}

\newcommand{\zgo}{\mathfrak{z}}

\newcommand{\Sgl}{\mathfrak{S}}

\newcommand{\N}{\mathbb{N}}
\newcommand{\Z}{\mathbb{Z}}
\newcommand{\A}{\mathbb{A}}
\newcommand{\C}{\mathbb{C}}

\newcommand{\Q}{\mathbb{Q}}
\newcommand{\R}{\mathbb{R}}

\newcommand{\Gm}{\mathbb{G}_m}
\newcommand{\Ga}{\mathbb{G}_a}

\newcommand{\calA}{\mathcal{A}}

\newcommand{\calE}{\mathcal{E}}
\newcommand{\calF}{\mathcal{F}}

\newcommand{\calH}{\mathcal{H}}

\newcommand{\calO}{\mathcal{O}}
\newcommand{\calP}{\mathcal{P}}

\newcommand{\calR}{\mathcal{R}}

\newcommand{\calT}{\mathcal{T}}

\newcommand{\nin}{\notin}

\newcommand{\tlphi}{\tilde \phi}

\newcommand{\tlGam}{\widetilde \Gamma}

\newcommand{\tlG}{\widetilde{G}}

\newcommand{\htau}{\widehat \tau}

\newcommand{\hDelta}{\widehat \Delta}

\newcommand{\brtau}{\bar \tau}
\newcommand{\brhtau}{\bar \htau}
\newcommand{\brago}{\bar \ago}
\newcommand{\bragop}{\overline{\ago_P^{+}}}
\newcommand{\bragoq}{\overline{\ago_Q^{+}}}
\newcommand{\bragor}{\overline{\ago_R^{+}}}
\newcommand{\brzgo}{\bar \zgo}
\newcommand{\brzgop}{\overline{\zgo_{P}^{+}}}
\newcommand{\brzgoq}{\overline{\zgo_{Q}^{+}}}

\newcommand{\sgn}{\mathrm{sgn}}

\newcommand{\rint}{\mathrm{rint}\,}

\newcommand{\treg}{\calT_{reg}}
\newcommand{\tzer}{\calT}

\newcommand{\tzers}{T_{G}}
\newcommand{\tregs}{T_{G,reg}}

\numberwithin{equation}{section}
\begin{document}
\title{Periods of automorphic forms over reductive subgroups}
\date{\today}
\author{Micha\l{} Zydor}
\address{Department of Mathematics\\ University of Michigan\\ Ann Arbor, MI USA}
\email{zydor@umich.edu}

\begin{abstract}
We present a regularization procedure of period integrals of automorphic forms 
on a group $G$
over an arbitrary reductive subgroup $G' \sbs G$. 
As a consequence we obtain an explicit $G'(\A)$-invariant functional 
on the space of automorphic forms on $G$ whose exponents avoid certain prescribed hyperplanes.
We also provide a necessary and sufficient condition for convergence of period integrals of automorphic forms 
in terms of their exponents.
\end{abstract}

\maketitle 
\section*{Introduction}\label{sec:intro}

\subsection{Our result}
Let $G$ be a connected reductive algebraic group and $G' \sbs G$ a subgroup, both defined over a number field $F$. 
An automorphic period is the integral
\[
\calP(\phi) := \int_{G'(F) \bsl G'(\A)}\phi(g)\,dg, \quad \phi \in \calA(G)
\]
where $\calA(G)$ is the space of automorphic forms on $G$. 

Periods appear in many areas related to number theory. They often pertain to special values of automorphic 
$L$-functions and can also be used to study automorphic representations. 

This article addresses the question of convergence of period integrals. 
It was proven in \cite{beuz2} that automorphic periods converge for cuspidal 
automorphic representations whenever $G/G'$ is a quasi-affine variety. 
However, unless the automorphic quotient $G'(F) \bsl G'(\A)$ is compact, 
the period integral is ill defined on the space of all automorphic forms $\calA(G)$. 

Nevetherless, it was first observed by Zagier \cite{zagier}, that in the realm of automorphic forms, 
divergent integrals can be regularized, providing a meaningful extension of the functional $\calP$, 
defined initially only for cusp forms, to the space $\calA(G)$. Zagier's work concerned the case 
of the scalar product on $GL_2$ and was further reinterpreted by Casselman \cite{cass}. 
However, it is the subsequent work of Jacquet-Lapid-Rogawski \cite{jlr} that provided 
 a framework that inspired our work. 
 
 We consider the case when $G'$ is a connected reductive subgroup. 
Let us describe briefly our approach. 
We define, what's commonly called a mixed truncation operator 
\begin{equation}\label{eq:mixTr}
\La^{T} : \calA(G) \to \text{ rapidly decreasing functions on }G'(F) \bsl G'(\A).
\end{equation}
The operator is called mixed as opposed to the standard one used by Langlands and Arthur \cite{arthur5}. 
Here, $T$ is a parameter in an $\R$-vector space of dimension equal to the split rank of $G'$. 
We obtain thus a family of operators on $\calA(G)$ that send it to the space of rapidly decreasing functions on 
$G'(F) \bsl G'(\A)$. 
Once such operator is defined and its properties established, cf. Section \ref{sec:relRedTh}, we can consider the 
integral
\[
\calP^{T}(\phi) := \int_{G'(F) \bsl G'(\A)}\La^{T}\phi(g)\,dg, \quad \phi \in \calA(G).
\]
The insight of \cite{jlr}, is that such an integral, even though not canonically defined, has a canonically 
defined constant term in the parameter $T$, at least for almost all automorphic forms. 

One defines then, the so called regularized period
\[
\calP(\phi) := \text{ constant term of } T \mapsto \calP^{T}(\phi).
\]
For $x \in G(\A)$ let $\phi_{x}(y) = \phi(yx^{-1})$. 
The main result of the paper can be stated then as follows. 

\btheo[cf. Theorem \ref{theo:theProp} ] 
For $\phi$ in a a properly defined subspace $\calA(G)^{*}$ of $\calA(G)$
we have
\[
\calP(\phi_{x}) = \calP(\phi), \quad \forall \ x \in G'(\A)
\]
and 
$\calP(\phi) $ is independent of any choices. 
\etheo

Moreover, Theorem \ref{thm:intCrit} provides a necessary and sufficient condition for period integrals of automorphic forms to converge. 
The article considers in fact a slightly more general setting where an automorphic 
character on the group $G'(\A)$ is included. In Paragraph \ref{ssec:nonConn} we briefly explain how to extend the 
regularization to a case where $G'$ is not connected.

\subsection{On the proofs}

The core of the proof lies in establishing the rapid decay property 
of $\La^{T}$. In essence, we follow the approach of Arthur \cite{arthur5}. 
We hope that the notation we employ makes this connection apparent. What allows us to follow Arthur's approach, 
and also the one of Jacquet-Lapid-Rogawski is a careful definition of the operator $\La^{T}$. 
Here's the definition
\begin{equation}\label{eq:mixTr2}
\La^{T}\phi (x) = \sum_{P \in \calF^{G}(P_{0}')}\varepsilon_{P}^{G}\sum_{\delta \in (P \cap G') \bsl G'(F)}
\htau_{P}(H_{0'}(\delta x)- T)\phi_{P}(\delta x).
\end{equation}
It should come as no surprise that $\La^{T}$ is defined as a certain alternating sum (we have $\varepsilon_{P}^{G} = \pm 1$)
of constant terms of $\phi$ ($\phi_{P}$ is the constant term of $\phi$ with respect to the parabolic subgroup $P$) 
truncated (via $\htau_{P}$) and summed to make everything $G'(F)$-invariant. 
The group $P_{0}'$ is a fixed minimal parabolic subgroup of $G'$. The set $\calF^{G}(P_{0}')$ is then 
defined as the set of parabolic subgroups of $G$ that can be defined via a cocharacter $\la : \Gm \to G$ 
that has values in $G'$ and regarded as a homomorphism $\Gm \to G'$ defines a parabolic subgroup of $G'$ 
containing $P_{0}'$. In particular, the definition ensures that $P \cap G' $ contains $P_{0}'$ 
for $P \in \calF^{G}(P_{0}')$. The function $H_{0'}$ is the Harish-Chandra function on $G'(\A)$ 
taking values in the space $\ago_{0'} = \Hom_{F}(\Gm, A_{0}') \otimes_{\Z} \R$ 
where $A_{0'} \sbs P_{0'}$ is a fixed $F$-split torus of $G'$. Finally, $\htau_{P}$ is the characteristic function of the dual cone 
of the cone $\ago_{0'} \cap \ago_{P}^{+}$, $\ago_{P}^{+}$ being the standard positive chamber associated 
to $P$ (and a fixed $F$-split torus of $P$ determined by $A_{0'}$).

The Section \ref{sec:cones} studies the various truncation functions, such as $\htau_{P}$ above. 
We took an approach of the theory of polyhedral cones. Let us cite \cite{barvinok} 
for an excellent reference to this theory. It turns out that this framework, albeit a little unorthodox in this setting, 
meets perfectly the needs of a regularization procedure. In the classical work of Arthur, there are two principal 
combinatorial constructions - these are the functions $\Gamma$ and $\sigma$ . Both are defined as certain alternating 
sums of characteristic functions of cones. In \cite{zhengZyd}, H. Zhang and the author study the generalization 
of the $\Gamma$ function. This work is recalled in Paragraph \ref{ssec:gamma}. 
The generalization of the $\sigma$-function is studied in Paragraph \ref{ssec:sigma}. 

Once the combinatorial properties of general truncation functions are established in Section \ref{sec:cones}, 
we start working in the setting of reductive groups. In Section \ref{sec:redTh} we set up the notation and study the reduction theory. 
We do not really prove anything new about the reduction theory of algebraic groups. We mostly recall some classical statements, 
at times prove minor improvements. 
It is in Section \ref{sec:relRedTh} that we properly introduce the relative setting $G' \sbs G$. 
We define the truncation operator $\La^{T}$ and prove some typical combinatorial properties it satisfies. 
The core result is the rapid decay property \eqref{eq:mixTr} of Theorem \ref{thm:rapDec}.
Before we can prove it however, we prove the "relative decomposition of $1$" result in Proposition \ref{prop:decOf1Rel}.
We provide some pictures of truncated cones throughout the text to provide some insight into our constructions. 
The final section \ref{sec:invPeriod} reaps the profits of the work performed in the former sections. We define 
the period as explained above and verify briefly their properties. 
Paragraph \ref{ssec:cuspEis} spells out our construction in the case of cuspidal Eisenstein series associated to maximal parabolic subgroups of $G$. Finally, in the last paragraph \ref{ssec:intCrit}, we prove Theorem \ref{thm:intCrit} 
which provides a definite criterion for convergence of periods of automorphic forms over reductive subgroups.. 

\subsection{Example: symmetric periods}

Let $\theta$ be an involution on $G$ and suppose $G'= G^{\theta}$. 
Fix $A_{0}$ a maximal $F$-split and $\theta$-stable torus of $G$ such that $A_{0}' = A_{0}^{\theta}$ is a 
maximal $F$-split torus of $G'$ (existence is proven in \cite{helWan}). 
We fix also $P_{0'}$ a minimal parabolic subgroup of $G'$ containing $A_{0}'$. 
We have then $\ago_{0'} = \ago_{0}^{\theta}$ and the set $\calP^{G}(P_{0}')$ in \eqref{eq:mixTr2} is the set of parabolic subgroups of $G$ 
that are $\theta$-stable and whose intersection with $G'$ contains $P_{0}'$.

\subsection{Applications and comparison with previous works}

We provide here a brief, and by no means exhaustive list of cases 
where regularized periods appeared in the literature (defined using a mixed truncation or not) 
and discuss some applications.

We've already mentioned the pioneering work of \cite{zagier, cass, jlr}. 
The paper \cite{jlr} was succeeded by \cite{lapRog} where general Galois periods 
are considered, as opposed to the case of $GL_n$ treated in \cite{jlr}. 
The regularized periods of \cite{lapRog} appear in the fine spectral expansion 
of Jacquet's relative trace formula established in  \cite{lapid} and finally, in the work of Feigon, Lapid and Offen \cite{feiLaOff}, 
find application to the study of 
periods of the form $U(n) \sbs GL$, where $U(n)$ is a unitary group.   
The mixed truncation operator of \cite{jlr} is also used by Yamana \cite{yamana2} 
to study (convergent) periods of residual automorphic forms.

Relative truncation has also been employed to study the periods of the form $Sp_{2n} \sbs GL_{2n}$
in the work of Offen and Yamana \cite{offen1, offen2, yamana1}. Note that some regularization process 
is indispensable to study this case as periods of cusp forms vanish by the work of Jacquet and Rallis 
\cite{jacqRall2}. 

Another place in the literature where the mixed truncation appears is in connection to the global Gan-Gross-Prasad conjecture. The work of Ichino-Yamana \cite{ichYamGl, ichYamUn} handle the case of the Rankin-Selberg integrals
 of type $GL_n \sbs GL_{n} \times GL_{n+1}$ and $U(n) \sbs U(n) \times U(n+1)$ respectively. 
 Similar constructions were used in author's work on the Jacquet-Rallis relative trace formula \cite{leMoi2, leMoi3, leMoi}.
 Our construction is closely related to theirs and can be proved to produce the same regularized periods.
The shapes of cones that are used is nonetheless different. 
This difference played a role in our work with Chaudouard \cite{chaudZyd} where the truncation 
studied here needed to be introduced in order to perform semi-simple descent of certain geometric distributions.
 
 Periods of automorphic forms have often relation to special values of automorphic $L$-functions. 
 In fact, many known cases pertain to special values of Langlands-Shahidi $L$-functions 
 and the special values are closely related to poles of intertwining operators. This fact can be used 
 by analyzing periods of specially chosen Eisenstein series. It is practically a feature 
 of such periods to be divergent, and a regularization of the integral reveals relation to special values. 
 This method has been employed by many authors  \cite{jiang, gjr1, polWanZyd}, 
 most notably Jiang, Ginzburg and Rallis, where Arthur's truncation operator was used. 
 The above mentioned work of Ichino-Yamana \cite{ichYamUn} also follows this path using a mixed truncation. 
 In a future collaboration \cite{polWanZyd2} with A. Pollack and C. Wan, we will 
apply the mixed truncation and regularization developed here to several new cases that emerged 
from considering the conjectures of Sakellaridis-Venkatesh \cite{sakVen}. The explicit formula 
of Corollary \ref{cor:cuspEis} is particularly useful for such questions. 

A different approach to divergent periods of automorphic forms is through analyzing periods of 
pseudo-Eisenstein series.
The work of Ginzburg-Lapid \cite{ginLap} studies the problem from this perspective. 
In  \cite{lapRog, ichYamGl} they are studied along regularized periods.
The recent work \cite{lapOff} uses them 
to formulate and address the problem of the $G'$-distinguished spectrum in $L^{2}(G(F) \bsl G(\A))$, 
focusing on $Sp(n) \times Sp(n) \sbs Sp(2n)$.

The regularized periods should naturally appear in the spectral development of a relative trace formula. 
As of now, the work of Lapid \cite{lapid} is the only reference in general rank, but 
regularized periods are also visible in numerous works on relative trace formula of Jacquet 
and collaborators \cite{jacq1, jacq2, jacqLai, jacqZag, jacqLaiRal}. 
The mixed truncation of Ichino-Yamana and Jacquet-Lapid-Rogawski appears also in 
our work on the Jacquet-Rallis relative trace formula \cite{leMoi3, jlr} and 
the corresponding regularized periods are expected to appear in the fine spectral expansion of this formula 
which is an ongoing work of P.-H. Chaudouard and the author. 
It should be noted that once a method of truncating periods is available, it suggests a way of 
truncating (and regularizing) the corresponding relative trace formula although this is something that 
needs more investigation and we do not address this in the present article. 

The question of defining the relative trace formula on a suitable Schwartz space 
was addressed in great generality in the work of Sakellaridis in \cite{sak1}. 
The second part of loc. cit. deals with the question of defining (or suggesting when it is not possible by geometric means)
 invariant distributions that appear in the geometric  
expansion of the relative trace formula in essentially the same generality as this article. 
In particular, the author defines regularized integrals of certain theta series 
through a novel regularization that he calls an evaluation map.
Even though the construction in \cite{sak1} does not use truncation, it does rely on 
the analytic continuation of integrals of exponential functions over polyhedral cones, just like our method does. 
Most notably though, the data that Sakellaridis uses in his regularization process, a certain fan of polyhedral cones, 
corresponds quite directly to the one we use in our construction (c.f \eqref{eq:bigFan}). 
This has been verified by Sakellaridis and the author in May 2018 during our stay at 
the Institute for Advanced Study. 

Chaudouard in \cite{chaud1} introduced a new way of truncating 
the notoriously difficult unipotent orbital integrals in the context of Arthur-Selberg trace formula 
and he obtained new expressions for the unipotent orbital integrals involving zeta functions and integrals of Eisenstein series.
The construction in \cite{chaud1} adapts the truncation of certain stacks of coherent sheaves constructed by 
Schiffman \cite{schiff}. On a more conceptual level, Chaudouard's approach 
is it to adapt methods of the very geometric theory of automorphic forms over function fields to the number field setting. 
For example, it is well known, that the partition of the Siegel domain proved by Arthur \cite{arthur3}, Lemma 6.4 
(whose generalization we prove in Proposition \ref{prop:decOf1Rel}) 
can be interpreted over function fields in terms of the Harder-Narasimhan filtration of vector bundles. 
It would be interesting to pursue this analogy in the relative setting but we restrict ourselves to the number field case. 
It should also be noted that the work of Levy \cite{levy} deals with a very general divergent integrals in the automorphic setting (infinitesimal trace formula for a group $G$ acting on a finite dimensional vector space), but the results in loc. cit. do not produce invariant distributions and it is not clear 
to us what is the precise relation of Levy's constructions to ours. 
 
Let us note finally that it is still an open question of how to best deal with invariant distributions on automorphic forms. As mentioned above, Sakellaridis \cite{sak1} has proposed a much more conceptual approach
with a view towards comparison of different relative trace formulae 
and attacking the beatiful conjectures of Sakellaridis and Venkatesh \cite{sakVen}. 
Michel and Venkatesh  \cite{miVen} perform yet another regularization procedure via convolutions with 
local measures in the context of $GL_2 \sbs GL_2^{3}$ and their approach is certainly prone to generalizations.
Lapid and his collaborators use pseudo-Eisenstein series circumventing the need for truncation as well. 
Kudla, Rallis \cite{kudRal} and Ichino \cite{ichino1} use, respectively, the enveloping and Hecke algebra 
action to obtain the regularized Siegel-Weil formula. 
We also mention an interesiting recent work of Wu \cite{wu} where 
a regularization (in the $PGL_2 \sbs PGL_2$ case) that allows for critical exponents is introduced. 
Even as 
truncation method goes, different approaches are possible. Arthur's truncation has been employed in the relative setting as explained above and it remains the only truncation tool for studying periods
over non-reductive subgroups. 
Quite remarkably, in \cite{lapRog}, the authors show that the Galois periods 
truncated using the mixed truncation are equal to those truncated using the standard (Arthur's) one.
The verdict will lie of course in the applications of these methods to the study of $L$-functions and 
the distinction problem.
It seems to us nonetheless that the proposed here truncation has the feature of being 
optimal, in the sense that it truncates only that what needs to be truncated and 
it does it in a clean way. 
For example, our analysis leads to establishing, in Theorem \ref{thm:intCrit}, 
a necessary and sufficient condition 
for convergence of (reductive) periods of automorphic forms.  
\subsection{Acknowledgments}
The author would like to thank  Pierre-Henri Chaudouard, Erez Lapid, Yiannis Sakellaridis and Chen Wan for 
discussions and useful remarks regarding this work. 
We thank Erez Lapid for suggesting that the 
convergence criterion for period integrals (Theorem \ref{thm:intCrit}) should be proved
along the same lines as the case of the scalar product considered in \cite{moeWald2}, Lemma I.4.11. 
This work began when the author was a member at the Institute for Advanced Study in 2018. 
A part of this work was also accomplished during our stay at the Institute for Mathematical Sciences at the National University of Singapore.
We would like to thank both of these institutions for their hospitality and pleasant working conditions.

\section{Polyhedral cones}\label{sec:cones}
This section studies polyhedral cones and certain 
functions that arise as signed sums of characteristic functions 
of such cones (in fact, we study two such functions, the $\Gamma$ function in \ref{ssec:gamma} 
and the $\sigma$ function in \ref{ssec:sigma}). The motivation stems of course from reduction theory that allows 
to reduce automorphic integrals to, essentially, integrals over cones. 
The tools we introduce mimic the constructions introduced by Langlands and Arthur.
In the context of the scalar product one is interested in cones associated to 
(semi) standard parabolic subgroups. The positive chambers are a very peculiar class of cones, in particular 
they are simplicial which means that the poset structure of the set of their faces is that of a 
set of subsets of a fixed set (the set of simple roots). 
In our approach to period integrals we will deal with more general cones. 
The cones we encounter are obtained as intersections of simplicial cones with subspaces. To the best of our knowledge 
such cones are essentially arbitrary. This lead us to extending the objects defined by Langlands and Arthur.

\subsection{The algebra of polyhedra}

Let $V$ be a finite dimensional Euclidean space over $\R$. Let $\bilif$ be the scalar product on it. 
By a half-space $\calH^{+}$ in $V$ we mean any subset of $V$ of the type 
\[
\calH^{+} = \{H \in V \ | \ \langle H,H \rangle \ge 0\}
\]
for some $H \in V \smin \{0\}$. 
A cone in $V$ is a finite intersection  of half-spaces of $V$.
Let $C$ be a cone. Then $-C = \{-H \ | \ H \in C\}$ is also a cone.

By a face of $C$ we mean an intersection of $C$ 
with a half-space containing $C$ or $-C$. 
We denote $\calF(C)$ the set of faces of $C$. 
Elements of $\calF(C)$ are also cones and the set $\calF(C)$
has a natural poset structure induced by inclusion relation.

In this section, $C$ will denote a cone in $V$, 
and we will use letters $F$, $E$ and $G$ to denote its faces. 
Note that $\calF(C)$ is closed under taking intersections. 
Let $F_{0} = F_{0}(C)$ be the minimal face of $C$. 
It is a linear subspace of $V$. 

Let $V_{C}$ be the subspace of $V$ composed of linear combinations of elements $C$. 
We have the following result, proved in Theorem 1.6 \cite{brunGub}.

\brop\label{prop:brunGub}
Let $C \sbs V$ be a cone such that $V_{C}= V$.
Let $\calH_1^{+}$, $\calH_2, \ldots, \calH_n^{+}$ 
be half spaces of $V$ 
such that 
\[
C = \bigcap_{i=1}^{n} \calH_{i}^{+}
\]
in an irredundant way. Then, the half-spaces $\calH_{i}^{+}$ are unique. 
\erop

We also note $V^{C}$ the orthogonal complement of $V_{C}$ in $V$. 
For $F \in \calF(C)$ we also note 
\[
V_{C}^{F} := V_{C} \cap V^{F}.
\]

Define
\[
\varepsilon_{C} = (-1)^{\dim V_{C}}
\]
and for $F \in \calF(C)$
\[
\varepsilon_{C}^{F} = (-1)^{\dim V_{C} - \dim V_{F}}.
\]
We have then the Euler characteristic formula
\begin{equation}\label{eq:euler}
\sum_{F \in \calF(C)} 
\varepsilon_{C}^{F} = 
\begin{cases}
1, & \text{ if }C\text{ is a subspace of }V, \\
0  & \text{ else}. 
\end{cases}
\end{equation}

The dual cone of the cone $C$ is defined as
\[
C^{\vee} =  \{ H \in V \ | \ \langle H , C \rangle \ge 0\}.
\]

\brem A cone is called non-degenerate if it has a non-empty 
interior in $V$ and if it contains no lines (i.e. its minimal face is $\{0\}$).
It is important for our construction to consider all cones in $V$, in particular the degenerate ones, 
if only to give a unified treatment to all faces of a non-degenerate cone.  
For example
\begin{itemize}
\item[-] A proper subspace $W \sbs V$ is an example of a cone. We have then that $W^{\vee}$ is the orthogonal complement of 
$W$ in $V$.
\item[-] 
Cones that are not pointed, i.e. whose minimal face is a non-trivial subspace of $V$. 
For example, if $C$ is such a cone and it is additionally of full dimension 
then the dual cone $C^{\vee}$ is a cone contained in a proper subspace of $V$. 
See picture \ref{figure:1}. 
\begin{flushleft}
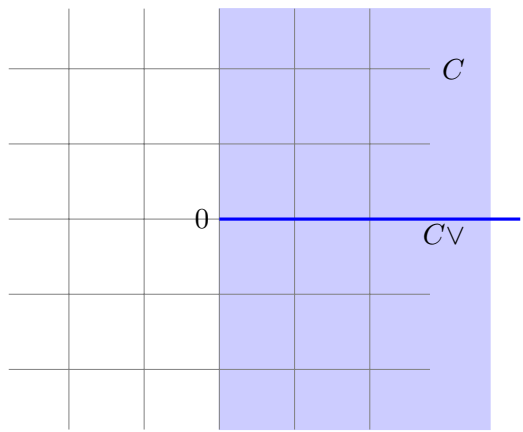
\begin{figure}
\begin{tikzpicture}[scale=2]
\filldraw[color=blue!20] (0,-1.4) rectangle (1.8,1.4);
\draw[step=.5cm,gray,very thin] (-1.4,-1.4) grid (1.4,1.4);
\draw[very thick, blue] (0,0) -- (2,0);
\draw (0,0) node[left] {$0$};
\draw (1.7,-0.1) node[left] {$C\vee$};
\draw (1.7,1) node[left] {$C$};
\end{tikzpicture}
  \caption{\tiny{Degenerate cone $C$ - a half space - and its dual $C^{\vee}$ - a half line}}
    \label{figure:1}
\end{figure}
\end{flushleft}
\item[-]  Cones $C$ contained in proper subspaces of $V$, that is such that $V_{C} \neq V$.
Note that since $(C^{\vee})^{\vee} = C$ this is a dual situation to the one discussed above and the figure \ref{figure:1} 
presents as well the relation of such cones to their duals. 
\end{itemize}
\erem 

Let $\rint C$ denote the relative interior of $C$. 
It is defined as the largest open subset of $V_{C}$ contained in $C$.
We will only apply it to closed cones, which should make the meaning of expressions as $\rint C^{\vee}$ clear: 
we have $\rint C^{\vee} = \rint (C^{\vee})$.

Let $F \in \calF(C)$. Define 
the angle cone as
\[
A(F,C) = \{ \la(H-Z)  \ | \ \la > 0, \, H \in C \} 
\]
where $Z$ in any point in $\rint F$ (the cone is independent of the choice of $Z \in \rint F$). See figure \ref{figure:2} 
for a typical example. 

\begin{flushleft}
\begin{figure}
\begin{tikzpicture}[scale=2]
\filldraw[color=blue!20] (2,2) -- (0,0) -- (2,0);
\draw[step=.5cm,gray, thin] (-1.4,-1.4) grid (1.7,1.7);
\draw[xstep=.1cm,ystep=3.cm, gray, very thin] (-1.5,0) grid (1.9,1.8);
\draw[rotate=225, xstep=.1cm,ystep=3.cm, gray,very thin] (-2.6,0) grid (1.6,1.4);
   \draw (0,0) node[left] {$0$};
   \draw (1.6,1.5) node[left] {$F_1$};
   \draw (1.7,0) node[left] {$F_2$};
   \draw (1.7,0.7) node[left] {$C$};
\draw (0.7,-1) node[left] {$A(F_1, C)$};
\draw (-0.4,1) node[left] {$A(F_2, C)$};
\end{tikzpicture}
  \caption{\tiny{Cone $C$, in blue, its two maximal faces $F_1, F_2$  and their angle cones given by half-spaces. }}
    \label{figure:2}
\end{figure}
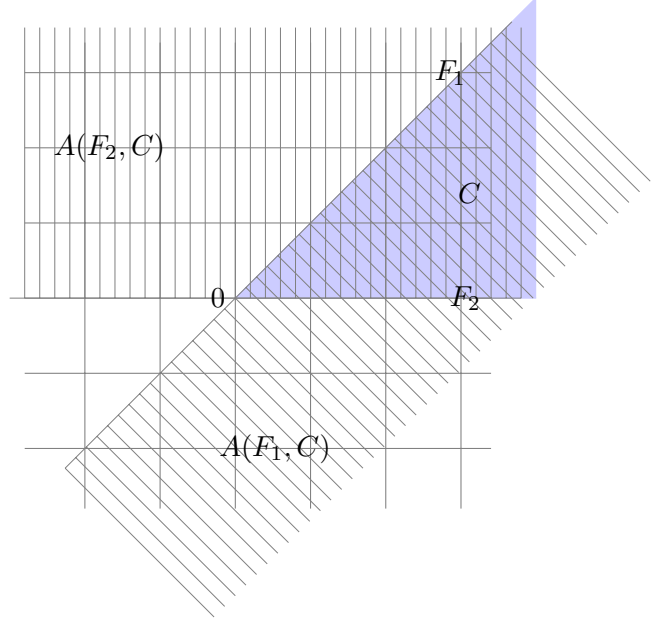
\end{flushleft}

We collect some useful facts about angle cones and their duals in 
the following lemma. 

\blem\label{lem:conProps} Let $C$ be a cone, $F_0$ the minimal face of $C$ and $F$ an arbitrary face of $C$. 
The following assertions hold. 
\begin{enumerate}
\item $A(F_{0},C) =C$ and $A(C,C) = V_{C}$.
\item There is an inclusion-preserving bijection
\begin{equation*}
	\begin{split}
	 \{\,G\in \calF(C): G\supset F \,\} & \xrightarrow{\sim} \calF(A(F,C))\\
	   G &\mapsto A(F,G).
	\end{split}
\end{equation*} 
Furthermore, for any face $G\in \calF(C)$ containing $F$, we have $A(A(F,G), A(F,C)) = A(G,C)$.
\item 
There is an inclusion-reversing bijection
\begin{equation*}
	\begin{split}
	 \{\,G\in \calF(C): G\supset F \,\} & \xrightarrow{\sim} \calF(A(F,C)^{\vee})\\
	   G &\mapsto A(G, C)^{\vee}.
	\end{split}
\end{equation*} 
Moreover, for any face $G \in \calF(C)$ containing $F$, we have $A(A(G,C)^{\vee}, A(F,C)^{\vee}) = A(F,G)^{\vee}$.
\end{enumerate}	
\elem

For a subset $A \sbs V$ we denote $[A]$ its characteristic function. 
We have the following classical result, see \cite{schneider} Theorem 3.1 
for a nice proof.

\btheo[The combinatorial 
Brianchon-Gram-Sommerville relation]\label{thm:bgs} For all $F \in \calF(C)$ we have
\[
\sum_{E \sps F}\varepsilon_{C}^{E} [\rint A(E,C)] = [- A(F,C)]
\]
and
\[
\sum_{F \in \calF(C)}\varepsilon_{F}^{F_{0}} [\rint F^{\vee}] = [-C^{\vee}].
\]
\etheo

 Let $H \in V$. 
Denote $H_{C}$ and $H^{C}$ the projections of $H$ onto 
$V_{C}$ and its orthogonal respectively. 
The following result, which can be viewed as a generalization of the Combinatorial Langlands Lemma, 
is proven in \cite{zhengZyd}, see also \cite{schneider}.

\brop\label{prop:langLemma} We have the following identities for all $H \in V$
\begin{enumerate}
\item
\[
\sum_{F \in \calF(C)} 
\varepsilon_{C}^{F} [\rint F](H_{F})[\rint A(F,C)^{\vee}](H^{F}) = 
\begin{cases}
1 , & \text{ if }C\text{ is a subspace of }V, \\
0   & \text{ else}. 
\end{cases}
\]
\item
\[
\sum_{F \in \calF(C)} 
\varepsilon_{C}^{F} [\rint A(F,C)](H)  [\rint F^{\vee}](H) = 
\begin{cases}
1, & \text{ if }C\text{ is a subspace of }V\text{ and }H=0, \\
0 & \text{ else}. 
\end{cases}
\]
\end{enumerate}
\erop

\begin{figure}
    \begin{minipage}{0.45\textwidth}
\begin{tikzpicture}[scale=2]
\filldraw[color=blue!10] (0.5, 1.5) -- (-1, 0) -- (1.5,0);
\filldraw[color=blue!40]  (-1,0) -- (0,1) -- (0.5,0.5) -- (0.5, 0) -- cycle;
\draw (0.5,0.5) node[right] {$T$};
\draw (-1,0) node[left] {$0$};
\draw (0.6,1) node[right] {$C$};
\draw (-0.5,0.25) node[right] {$ \Gamma(C, \cdot ,T )$};
\draw[step=.5cm,gray, thin] (-1.6,-0.6) grid (1.2,1.5);
\end{tikzpicture}
  \caption{\tiny{$\Gamma(C, \cdot,T )$ is the characteristic  function of the dark blue region 
  seen as a subset of the interior of the cone $C$.}}
    \label{figure:3}
        \end{minipage} \hfill
        \begin{minipage}{0.45\textwidth}
\begin{tikzpicture}[scale=2]
\filldraw[color=blue!20] (0,-1.4) rectangle (1.8,1.4);
\draw[step=.5cm,gray,very thin] (-1.4,-1.4) grid (1.4,1.4);
\draw (0,0) node[left] {$0$};
\draw (1.7,-0.1) node[left] {$C$};
\filldraw[blue] (1,0.8) circle (0.03cm);
\draw (1.3,0.8) node[left] {$T$};
\draw (0.1,.95) node[right] {$ \Gamma(C, \cdot ,T )$};
\draw[thick, blue] (0,0.8) -- (1,0.8);
\end{tikzpicture}
  \caption{\tiny{For $C$ a half plane and $T \in C$, the function $\Gamma(C, \cdot ,T )$ is the characteristic 
  function of the blue interval (it is open on the left and closed on the right)}}
    \label{figure:4}
    \end{minipage}
\end{figure}

\subsection{The $\Gamma$ function}\label{ssec:gamma}

For $T \in V$, let
\[
\Gamma(C, H, T) = \sum_{F \in \calF(C)}\varepsilon_{F}^{F_{0}}
[\rint A(F,C)](H)[\rint F^{\vee}](H-T), \quad H \in V.
\]
This is a generalization to arbitrary cones of the $\Gamma$ function introduced by Arthur \cite{arthur2}.
The function should be seen as a function of $H$ and $T$ as the truncation parameter. 
The figure \ref{figure:3} presents the classical picture of the support of $\Gamma(C, \cdot , T)$ 
in the case of two-dimensional acute cone and $T \in C$. In this case it should be thought of as a truncated fundamental domain 
of a rank 2 reductive algebraic group. 
Arthur proved that such a picture generalizes to simplicial cones obtained as chambers of positive weights 
associated to parabolic subgroups. 
The figure \ref{figure:4} presents the $\Gamma$ function of a degenerate cone.
For arbitrary cones the $\Gamma$ 
function is usually not a characteristic function. 
Nevertheless, it is always compactly supported and exhibits some nice properties. 
We study it in our work  with H. Zheng \cite{zhengZyd}, we will limit ourselves here to stating the results 
that are needed for our construction.  

We start with the following observation

\blem\label{lem:gamma0} $\Gamma(C,H,T) = 0$ unless $H_{F_0} = T_{F_0}$ and $H \in V_{C}$.
\elem
\bdem
Let $F \in \calF(C)$. 
If $[\rint F^{\vee}](H - T) \neq 0$ we must have $H_{F_0} = T_{F_0}$ 
since $F_0$ is a vector space contained in $F$. On the other hand, we have $A(F,C) \sbs V_{C}$. 
\edem

The following lemma is proved in \cite{zhengZyd}, Lemma 2.5 and Proposition 2.9

\blem\label{lem:gamma1} Let $T \in V$.
\begin{enumerate}
\item For all $H \in V$, $\Gamma(C, H, T) = 0$ 
unless $\langle H, H-T \rangle \le 0$. 
\item The function $H \in V \mapsto \Gamma(C, H ,T)$ is compactly supported for all $T$. 
\end{enumerate}
\elem 

The next result can be seen as a geometric justification for introducing the $\Gamma$-function. It 
arises naturally when one looks at the difference of a cone and its translate by a fixed vector $T$. 
The figure \ref{figure:5} presents geometrically its content 
for a non-degenerate cone and the figure \ref{figure:5prim} for a generate one.  

\begin{figure}
 \captionsetup{width=.5\linewidth}
    \begin{minipage}{0.45\textwidth}
\begin{tikzpicture}[scale=2]
\filldraw[color=blue!20]  (1.5, 2.5) -- (-1, 0) -- (2,0);
\filldraw[color=blue!10]  (1.8, 1.8) -- (0.5, 0.5) -- (2, 0.5);
\filldraw[color=blue!40] (-1,0) -- (0,1) -- (0.5,0.5) -- (0.5, 0) -- cycle;
\draw (-1,0) node[left] {$0$};
\draw (0.2,1.2) node[left] {$F_1$};
\draw (1.2,0) node[below] {$F_2$};
\draw (0.5,0.5) node[right] {$T$};
\draw (-0.5,0.25) node[right] {$\Gamma(C, \cdot ,T)$};
\draw (1,0.8) node[right] {$C+T$};
\draw (0.8,0.3) node[right] {$A_2$};
\draw (0.8,1.4) node[right] {$A_1$};
\draw[step=.5cm,gray,very thin] (-1.4,-.4) grid (1.4,2.4);
\end{tikzpicture}
  \caption{\tiny{Decomposition of a an interior of a cone $C$ into its translate by $T$, the $\Gamma(C, H,T) \neq 0$ region, 
  and two regions $A_1$, $A_2$ whose respective indicators are 
  $\Gamma(A(F_1,C), H^{F_1}, T^{F_1})[\rint F_1](H_{F_1}-T_{F_1})$ 
  and 
  $\Gamma(A(F_2,C), H^{F_2}, T^{F_2})[\rint F_2](H_{F_2}-T_{F_2})$}
  }
    \label{figure:5}
        \end{minipage} \hfill
        \begin{minipage}{0.45\textwidth}
\begin{tikzpicture}[scale=2]
\filldraw[color=blue!40] (0,-1.4) rectangle (1,1.4);
\filldraw[color=blue!20] (1,-1.4) rectangle (2.4,1.4);
\draw[step=.5cm,gray,very thin] (-0.4,-1.4) grid (2.4,1.4);
\draw (0,0) node[left] {$0$};
\draw (2.1,.6) node[left] {$C + T$};
\filldraw (1,0.7) circle (0.03cm);
\draw (1.25,0.7) node[left] {$T$};
\draw (0,.35) node[right] {$ \Gamma(C, H^{F} ,T^{F} )$};
\end{tikzpicture}
  \caption{\tiny{Decomposition of the interior of a half plane $C$ into its translate by $T$ 
and a horizontally bounded region $\Gamma(C, H^{F},T^{F} ) \neq 0$. }}
    \label{figure:5prim}
    \end{minipage}
\end{figure}

\blem\label{lem:gamma2} We have the following identity
\[
[\rint C](H) = 
\sum_{F \in \calF(C)}
\Gamma(A(F,C), H^{F}, T^{F})
[\rint F](H_{F}-T_{F}), \quad H,T \in V.
\]
\elem
\bdem
We have
\begin{multline*}
\sum_{F \in \calF(C)}
\Gamma(A(F,C), H^{F}, T^{F})
[\rint F](H_{F}-T_{F}) = \\
\sum_{E \sps F}\varepsilon_{E}^{F}[\rint A(E,C)](H^{F})
[\rint A(F,E)^{\vee}](H^{F} - T^{F})[\rint F](H_{F}-T_{F}) = \\
\sum_{E}[\rint A(E,C)](H^{F})
\dsl 
\sum_{F \in \calF(E)}
\varepsilon_{E}^{F}
[\rint A(F,E)^{\vee}](H^{F} - T^{F})[\rint F](H_{F}-T_{F})
\rb = [\rint C](H) 
\end{multline*}
where we use part 2 of Proposition \ref{prop:langLemma} to obtain the last equality. 
\edem

Let us also note, an essentially formal, dual version of the Lemma above

\bcor\label{cor:gamma2} The following identity holds
\[
[\rint C^{\vee}](H-T) = 
\sum_{F \in \calF(C)}
\varepsilon_{F}^{F_{0}} [\rint A(F,C)^{\vee}](H^{F}) \Gamma(F, H_{F},T_{F}), \quad H,T \in V.
\]
\ecor
\bdem 
Note first that we have
\begin{equation}\label{eq:gamma1}
\Gamma(C,H,T) = \varepsilon_{C}^{F_{0}}\Gamma(C^{\vee}, H-T,-T).
\end{equation}
Indeed, the third point in Lemma \ref{lem:conProps}, implies that
\[
\varepsilon_{F}^{F_{0}}
[\rint A(F,C)](H)[\rint F^{\vee}](H-T) = 
\varepsilon_{C}^{F_{0}}\varepsilon_{A(F,C)^{\vee}}^{A(C,C)^{\vee}}
[\rint A(A(F,C)^{\vee}, C^{\vee})](H-T)[\rint (A(F,C)^{\vee})^{\vee}](H)
\]
which implies \eqref{eq:gamma1} by the definition of $\Gamma$. 

To prove the Corollary, let us apply Lemma \ref{lem:gamma2} to the cone $C^{\vee}$ setting $H = H-T$ 
and $T = -T$ therein. 
We obtain
\[
[\rint C^{\vee}](H-T) = 
\sum_{F \in \calF(C)}
\Gamma(A(A(F,C)^{\vee},C^{\vee}), (H-T)^{A(F,C)^{\vee}}, -T^{A(F,C)^{\vee}})
[\rint A(F,C)^{\vee}](H_{A(F,C)^{\vee}}).
\]
The result follows now from part 3 of Lemma \ref{lem:conProps}, the identity \eqref{eq:gamma1} 
and the fact that
\[
V_{A(F,C)^{\vee}} = V^{F}.
\]
\edem

The results mentioned so far should be thought of as classical, or absolute. The neHt result concerns the relative case. 
It is a refinement of the decomposition in Lemma \ref{lem:gamma2}. In addition to a cone $C$ and a truncation parameter $T$ 
we also have a subdivision of the cone $C$ into a fan of subcones. 

\blem\label{lem:gamma3} Let $C_1, \ldots, C_n$ be cones in $V$ 
such that $C = \bigcup_{i}C_i$ 
and such that $C_i \cap C_j$ is a face of $C_i$ and $C_j$ 
for all $i,j = 1,\ldots, n$. 
Then
\[
[\rint C](H) = 
\sum_{F \in \calF(C)}
\sum_{
\begin{subarray}{c}
G \in \bigcup_{i}\calF(C_i)\\
 \rint G \sbs \rint F 
\end{subarray}
}
\Gamma(A(F,C), H^{G}, T^{G})
[\rint G](H_{G}-T_{G}), \quad H,T \in V.
\]
\elem

This arises in the context of automorphic periods as follows. 
We have an inclusion of reductive algebraic groups $G' \sbs G$. We can fix maximal split tori of the respective groups 
so that we have the inclusion $\ago_{0'} \sbs \ago_{0}$  (see Paragraph \ref{ssec:chambs} for an unexplained notation). 
It is the space $\ago_{0'}$ that serves the role of the ambient space $V$ and the role of the cone $C$ 
is played by the cone $\ago_{0'}^{+}$ associated to a minimal parabolic subgroup of $G'$. 
Since $\ago_{0} = \bigcup_{P} \ago_{P}^{+}$ where the union runs over semi-standard parabolic subgroups of $G$ 
we obtain a decomposition of the cone $\ago_{0'}^{+}$ by taking intersections with this decomposition. 
This is exactly the fan of cones that we obtain. Interestingly, as mentioned in the introduction, 
a closely related fan is used in \cite{sak1}. 
The figure \ref{figure:6} depicts the content of Lemma \ref{lem:gamma3} geometrically in the case of a degenerate 
cone. It is important why this is pertinent to our case. In the setting of the scalar product, the center 
of the group is automatically discarded. That is why one works with non-degenerate cones in this setting. 
In the relative setting $G' \sbs G$, one can only mod out by the intersection of the center of $G$ 
with the group $G'$. This means that sometimes the center of $G'$ plays a role. 
As an extreme, but non-trivial, case one can look at the case when $G'$ is a torus. 
For example when $GL_1 \sbs GL_2$ ($GL_1$ embedded into the upper left corner), 
$\ago_{0'}$ is one dimensional and 
the two semi-standard Borel subgroups 
of $GL_2$ partition it into a "positive" and "negative" part. 
Figure \ref{figure:6} presents this case.
In case of the diagonal inclusion $GL_2 \sbs GL_3$, the space $\ago_{0'}$ is two dimensional and the 
the positive chamber of $GL_2$ is a half space in it. Three chambers in $GL_3$ intersect this chamber and figure 
\ref{figure:7} presents the decomposition of Lemma \ref{lem:gamma3} in this case. 

\begin{flushleft}
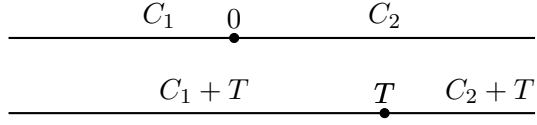
\begin{figure}
\begin{tikzpicture}[scale=2]
\draw[thick]  (-1.5, 0.5) -- (2, 0.5);
\filldraw (0,0.5) circle (0.03cm);
\draw (1,0) node[above] {$T$};
\draw (0,0.5) node[above] {$0$};
\draw[thick]  (-1.5, 0) -- (2, 0);
\filldraw (1,0) circle (0.03cm);
\draw (1,0) node[above] {$T$};
\draw (-0.5,0.5) node[above] {$C_1$};
\draw (-0.2,0) node[above] {$C_1 + T$};
\draw (1,0.5) node[above] {$C_2$};
\draw (1.7,0) node[above] {$C_2 + T$};
\end{tikzpicture}
  \caption{\tiny{The decomposition of a line into two half-lines and below it
 the content of Lemma \ref{lem:gamma3} in this context. The line 
 is decomposed into three disjoint regions, two open half-lines ending in $T$ and the singleton $\{T\}$. }
  }
    \label{figure:6}
\end{figure}
\end{flushleft}
\begin{flushleft}
\begin{figure}
\begin{tikzpicture}[scale=2]
\filldraw[color=blue!1] (0,1.5) -- (0,0) -- (1.5,1.5);
\filldraw[color=blue!15]   (1.5,-1.5) -- (0,0) -- (1.5,1.5);
\filldraw[color=blue!10]  (1.5,-1.5) -- (0,0) -- (0,-1.5);
\draw[step=.5cm,gray,very thin] (-0.4,-1.4) grid (1.4,1.4);
\draw (0,0) node[left] {$0$};
\draw (0.6,1) node[left] {$C_1$};
\draw (1.2,0) node[left] {$C_2$};
\draw (0.6,-1) node[left] {$C_3$};
\filldraw[color=blue!30]   (4,1.5) -- (4,0.3) -- (5,0.3) -- (5,1.5);
\filldraw[color=blue!45]   (4,-1.5) -- (4,0.3) -- (5,0.3) -- (5,-1.5);
\filldraw[color=blue!1] (5,1.5) -- (5,0.3) -- (6.5,1.5);
\filldraw[color=blue!15]   (6.5,-1.5) -- (5,0.3) -- (6.5,1.5);
\filldraw[color=blue!10]  (6.5,-1.5) -- (5,0.3) -- (5,-1.5);
\draw[step=.5cm,gray,very thin] (3.6,-1.4) grid (6.4,1.4);
\draw (4,0) node[left] {$0$};
\draw (5.8,1.1) node[left] {$C_1 + T$};
\draw (6.2,0.3) node[left] {$C_2 + T$};
\draw (5.8,-0.7) node[left] {$C_3 + T$};
\draw (5.3,0.3) node[left] {$T$};
\draw (4.9,0.42) node[left] {$\Gamma(C, \cdot, T)$};
\filldraw[blue] (5,0.3) circle (0.03cm);
\draw[thick, blue] (4,0.3) -- (5,0.3);
\end{tikzpicture}
  \caption{\tiny{On the left, the half space $C$ is decomposed into a union of three non-degenerate cones. On the right 
  the induced decomposition of $C$ of Lemma \ref{lem:gamma3} is shown. 
  The upper rectangular parts corresponds to the face of $C_1$ that intersects the minimal face of $C$ (i.e. the $y$-axis )
  and the lower rectangular region corresponds to the analogous face of $C_3$. The blue interval is precisely 
  $\Gamma(C, \cdot, T)$ - the term corresponding to the face $\{0\}$.} }
    \label{figure:7}
\end{figure}
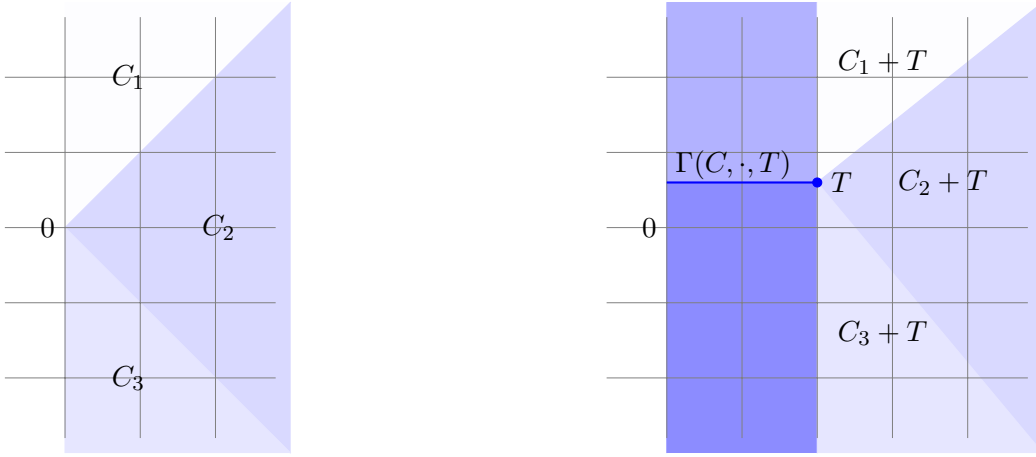
\end{flushleft}

\bdem
Using Lemma \ref{lem:gamma2} we've just proven, 
it is enough to prove that for all $F \in \calF(C)$ 
we have
\begin{equation}\label{eq:gammaOldNew}
\Gamma(A(F,C), H^{F}, T^{F})
[\rint F](H_{F}-T_{F}) = 
\ssum{
G \in \bigcup_{i}\calF(C_i)\\
 \rint G \sbs \rint F 
}\Gamma(A(F,C), H^{G}, T^{G})
[\rint G](H_{G}-T_{G}).
\end{equation}
Using the definition of $\Gamma(A(F,C), \cdot, \cdot)$ 
and the properties of angle cones,  
it is enough to prove for all $E \in \calF(C)$ 
containing $F$
the equality
\begin{multline*}
[ \rint A(E,C)](H^{F})[\rint A(F, E)^{\vee}](H^{F} - T^{F}) [\rint F](H_{F}-T_{F}) = 
 \\
\ssum{
G \in \bigcup_{i}\calF(C_{i}) \\
 \rint G \sbs \rint F 
}
[\rint A(E,C)](H^{G})[\rint A(F,E)^{\vee}](H^{G} - T^{G})  [\rint G](H_{F}-T_{F}) 
\end{multline*}

First of all, we have 
\[
[\rint A(E,C)](H^{F}) = [\rint A(E,C)](H) = [\rint A(E,C)](H^{G}) 
\]
since $G$ and $F$ are contained in $E$. 

Let $Z = H-T$. 
We have
\[
\rint F = \bigsqcup_{
\begin{subarray}{c}
G \in \bigcup_{i}\calF(C_i)\\
 \rint G \sbs \rint F 
\end{subarray}
} \rint G
\]
so 
\[
[\rint F](Z_{F})  = \sum_{
\begin{subarray}{c}
G \in \bigcup_{i}\calF(C_i)\\
 \rint G \sbs \rint F 
\end{subarray}
} [\rint G](Z_{F}).
\]

It remains to observe that the following equality holds
\[
[\rint A(F,E)^{\vee}](Z^{F}) [\rint G](Z_{F}) = 
[\rint A(F,E)^{\vee}](Z^{G}) [\rint G](Z_{G}).
\]
\edem

\subsection{The $\sigma$ function}\label{ssec:sigma}

In this section we study the generalization 
of Arthur's $\sigma_{P}^{Q}$ 
function, introduced in \cite{arthur3}, Section 6. It plays an important role in the proof of the rapid decay 
property of the truncation operator (Theorem \ref{thm:rapDec}).

For $F \in \calF(C)$, we define
\[
\sigma(F, C) = \sum_{E \sbs F} \varepsilon_{F}^{E}   [\rint A(E,C)]  [\rint E^{\vee}].
\]

We have the following two Lemmas that are easily established. 

\blem\label{lem:sigZero} 
We have
\[
\sigma(C, C)(H) = 
\begin{cases}
1, & \text{ if }C\text{ is a subspace of }V\text{ and }H=0, \\
0 & \text{ else}. 
\end{cases}
\]
\elem
\bdem 
The lemma follows from Proposition \ref{prop:langLemma}.
\edem

\blem\label{lem:htauTauSigSum} For all $F \in \calF(C)$ we have
\[
  [\rint A(F,C)]  [\rint F^{\vee}] = \sum_{E \sbs F}\sigma(E,C).
\]
\elem
\bdem 
We have by definition
\[
\sum_{E \sbs F}\sigma(E,C) = 
\sum_{G \sbs E \sbs F}
\varepsilon_{E}^{G}   [\rint A(G, C)]  [\rint G^{\vee}] = 
 \sum_{G \sbs F} [\rint A(G,C)]  [\rint G^{\vee}]
\sum_{G \sbs E \sbs F}\varepsilon_{E}^{G}. 
\]
Invoking the identity \eqref{eq:euler} completes the proof. 
\edem

For example, let's look at a non-degenerate cone $C$ in a two dimensional space. 
Let's take $F = C$ - in the lemma above. 
Taking Lemma \ref{lem:sigZero} into consideration we obtain
\[
[\rint C^{\vee}] = [\rint C] + \sigma(F_1, C) + \sigma(F_2, C)
\]
where $F_1, F_2$ are the maximal faces of $C$. The figure \ref{figure:8} 
depicts this decomposition in the case when $C$ is acute.

\begin{flushleft}
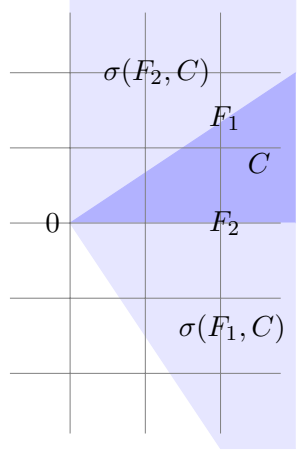
\begin{figure}
\begin{tikzpicture}[scale=2]
\filldraw[color=blue!10] (0,1.5) -- (0,0) -- (1.5,1) -- (1.5,1.5);
\filldraw[color=blue!30] (1.5, 0) -- (0,0) -- (1.5,1);
\filldraw[color=blue!10] (1.5, -1.5) -- (1,-1.5) -- (0,0) -- (1.5,0);
\draw[step=.5cm,gray,very thin] (-0.4,-1.4) grid (1.4,1.4);
\draw (0,0) node[left] {$0$};
\draw (1.2,0.7) node[left] {$F_1$};
\draw (1.2,0) node[left] {$F_2$};
\draw (1.4,0.4) node[left] {$C$};
\draw (1.5,-0.7) node[left] {$\sigma(F_1, C)$};
\draw (1,1) node[left] {$\sigma(F_2, C)$};
\end{tikzpicture}
  \caption{\tiny{The decomposition of the cone $C^{\vee}$ into the $C$ and two regions 
  $\sigma(F_1, C)$ and $\sigma(F_2, C)$ associated to its maximal faces $F_1$ and $F_2$.}}
    \label{figure:8}
\end{figure}
\end{flushleft}

The following is the main result concerning the $\sigma$ function. 
It's classical analogue is Lemma 6.1 and Corollary 6.2 of \cite{arthur3}. Our proof 
differs however of the one in loc cit. even in that case. 

\brop\label{prop:sigma} 
We have for $H \in V$
\[
\sigma(F, C)(H) \neq 0 \Longrightarrow [\rint A(F,C)](H) \neq 0.
\]
Moreover, there exist a constant $k > 0$ such that 
for all $H \in V$  such that $\sigma(F, C)(H) \neq 0$ 
we have
\[
\|H_{F}\| \le k\|H^{F}\|.
\]
\erop
\bdem 
The first part of the proposition is  clear since $A(E,C) \sbs A(F,C)$
for $E \in \calF(F)$. 

Let's prove the second assertion. We will prove it by induction on $|\calF(C)|$ - the cardinality of 
$\calF(C)$. If $|\calF(C)| = 1$, then $C$ is a vector space 
and $\sigma(C,C)(H) = 0$ unless $H = 0$ by Lemma \ref{lem:sigZero}
so the result follows. Assume the property
is proven for all cones $C$ 
such that $|\calF(C)| < n$, where $n \ge 2$.
Assume that $C$ is fixed and such that $|\calF(C)| = n$

Note first that if $H_{F} = 0$ the desired inequality is true for any positive constant $k$. 
It is enough thus to prove the claim when $H_{F} \neq 0$, which we will assume form now on.
Using Lemma \ref{lem:gamma2}
with $T = H^{F}$ we obtain
for all $E \in \calF(F)$
\[
  [\rint A(E,C)](H) = 
   \sum_{G \sps E}\Gamma(A(G,C), H^{G}, (H^{F})^{G})[\rint A(E,G)](H_{G} - (H^{F})_{G}).
\]
We observe that for $E \sbs F \cap G$ we have 
\[
[\rint E^{\vee}](H) = [\rint E^{\vee}](H_{G} - (H^{F})_{G}).
\]
Therefore we can rewrite $\sigma(F, C)(H)$ as
\begin{multline*}
\sum_{E \sbs F} \varepsilon_{F}^{E}   [\rint A(E,C)](H)  [\rint E^{\vee}](H) = \\
\sum_{E \sbs F}  \varepsilon_{F}^{E} \sum_{G \sps E}\Gamma(A(G,C), H^{G}, (H^{F})^{G})[\rint A(E,G)](H_{G} - (H^{F})_{G})
[\rint E^{\vee}](H_{G} - (H^{F})_{G})   = \\
\sum_{G \in \calF(C)} \varepsilon_{G}\Gamma(A(G,C), H^{G}, (H^{F})^{G})
\sigma(F \cap G, G, H_{G} - (H^{F})_{G}).
\end{multline*}
Fix $G \in \calF(C)$. 
It is enough to prove that 
\begin{equation*}
\Gamma(A(G,C), H^{G}, (H^{F})^{G})
\sigma(F \cap G, G, H_{G} - (H^{F})_{G}) \neq 0
\end{equation*}
implies the desired inequality. 

If $G = C$, we have 
\[
\sigma(F \cap G, G, H_{G} - (H^{F})_{G}) = 
\sigma(F,C)(H_{F}) = \sigma(F,F)(H_{F}) 
\]
which is zero by the assumption $H_{F} \neq 0 $ and Lemma \ref{lem:sigZero}.

Suppose $G \neq C$. 
Using $\sigma(F \cap G, G, H_{G} - (H^{F})_{G}) \neq 0$, 
the fact that $H = H_F + H^{F}$ and the induction hypothesis we have
$H_{G} - (H^{F})_{G} = (H_{F})_{G}$ and
\begin{equation}\label{eq:sig1}
\|(H_{F})_{F \cap G}\| \le k_0 \|(H_{F})^{F \cap G}_{G}\|
\end{equation}
for some $k_0 > 0$. If $F \sbs G$ this implies already the inequality we want. 

Assume then that $F \nsbs G$. 
Adding $ \|(H_{F})^{F \cap G}_{G}\|$ to the both sides of the inequality \eqref{eq:sig1} we obtain
\begin{equation}\label{eq:sig2}
\|(H_{F})_{G}\| \le k_1 \|(H_{F})^{F \cap G}_{G}\|
\end{equation}
for some constant $k_1 > 0$.

Since $\Gamma(A(G,C), H^{G}, (H^{F})^{G}) \neq 0$, 
we have, using Lemma \ref{lem:gamma1}
\[
\langle H^{G}, H^{G} - (H^{F})^{G}) \rangle \le 0.
\]
Using $H = H_F + H^{F}$ and the Cauchy-Schwarz inequality we obtain 
\[
\|(H_{F})^{G}\|^{2} \le - \langle (H^{F})^{G}, (H_{F})^{G} \rangle \le \|(H^{F})^{G}\| \|(H_{F})^{G}\|.
\]
Note that $(H_{F})^{G} = 0$ implies $H_{F} \in V_{F \cap G}$. 
This in turn means that the right hand member of inequality \eqref{eq:sig2} is zero 
which gives $(H_{F})_{G} = 0$ by this very inequality and consequently $H_{F} = 0$ 
which we have already excluded from considerations. It is legitimate therefore to deduce from the above inequality 
the inequality
\[
\|(H_{F})^{G}\| \le \|(H^{F})^{G}\|
\]
which implies that
\begin{equation}\label{eq:sig3}
\|(H_{F})^{G}\| \le \|H^{F}\|.
\end{equation}

We have 
$(H_{F})^{G} = (H_{F}^{F \cap G})^{G}$ and the orthogonal projection 
$V \to V^{G}$ restricted to $V_{F}^{F \cap G}$ 
is injective. By Lipschitz continuity of linear maps, this implies that there exists a constant $k_2$ such that
\begin{equation}\label{eq:sig4}
\|(H_{F})^{F \cap G}_{G} = (H_{F}^{F \cap G})_{G} \| \le \| H_{F}^{F \cap G} \| \le k_{2} \|(H_{F}^{F \cap G})^{G}  = (H_{F})^{G}\|.
\end{equation}

Applying \eqref{eq:sig4} to the right hand side of \eqref{eq:sig2} 
we obtain 
\begin{equation}\label{eq:sig5}
\|(H_{F})_{G}\| \le k_1 \|(H_{F})^{F \cap G}_{G}\| \le k_{1} k_2 \|(H_{F}^{F \cap G})^{G}  = (H_{F})^{G}\| \le 
k_1 k_2 \|H^{F}\|
\end{equation}
where we employed \eqref{eq:sig3} to get the last inequality. 

Finally, adding the inequalities \eqref{eq:sig3} and \eqref{eq:sig5} 
we obtain 
\[
 \|H_{F}\| = \|(H_{F})_{G}\| + \|(H_{F})^{G}\| \le k_1 k_2 \|H^{F}\| + \|H^{F}\| = (k_1 k_2 + 1)\|H^{F}\|
\]
as desired.
\edem

\subsection{Fourier transform}\label{ssec:FT}

 We collect here some facts about Fourier transforms of cones, and more general polyhedra. 
The proofs can be found in \cite{lawrence} and \cite{khovPukh} 
as well as in the beautiful book \cite{barvinok}.

We fix the unique Haar measure on $V$ and all its subspaces, giving the volume $1$ to any fundamental parallelepiped. 
We also note $V_{\C} := V \otimes_{\R}\C$ and extend the bilinear product $\bilif$ to it in a natural way.
For $\la \in V_{\C}$ we denote by $\Rel(\la)$ the real part of $\la$.

Let $\calF^{min}(C)$ denote the subset of $\calF(C)$ 
consisting of faces of dimension $\dim F_{0} + 1$.

Let $q \in \C[V]$ be a polynomial on $V$. Consider the integral
\[
\scrF(C, q,\la) := 
\int_{V_{C}^{F_{0}(C)}}[C](H)e^{\langle \la, H \rangle}q(H)\,dH,
 \quad \la \in V_{\C}.
\]
The integral defining $\scrF(C, q,\la)$ converges absolutely for
\[
\Rel(\la) \in -\rint (V^{F_{0}} \cap C)^{\vee}
\]
and admits a meromorphic continuation to $V_{\C}$, denoted by $\scrF(C, q,\la)$ as well. 
Moreover, the function $\scrF(C, q,\la)$ is holomorphic 
outside of the closed set of $\la \in V_{\C}$ vanishing on one of the lines $V_{F}^{F_{0}}$ 
for $F \in \calF^{min}(C)$.

For further reference, we define the open set $V_{\C}^{C-reg}$ of $V_{\C}$ as follows 
\[
V_{\C}^{C-reg} := \{ \la \in V_{\C} \ | \ \langle \la, V_{F}^{F_{0}} \rangle \neq 0 \ \forall F \in \calF^{min}(C)\}.
\]
We have then that $\scrF(C, q,\la)$ is holomorphic for $\la \nin V_{\C}^{C-reg}$.

A polynomial exponential function on $V$ is a function $f$ of the form
\[
f(T) = \sum_{\la \in V_{\C}}e^{\langle \la, T \rangle }q_{\la}(T)
\]
where $q_{\la} \in \C[V]$ with $q_{\la} = 0$ for all but finitely many $\la \in V_{\C}$.
Such a decomposition is unique.
The purely polynomial part of $f$ is by definition the polynomial $q_{0}$ 
corresponding to $\la = 0$.

In the notation above, for $T \in V$ let 
\begin{equation}\label{eq:FTGamma}
\scrF(\Gamma(C), T,q,\la) := \int_{V_{C}^{F_{0}(C)}}\Gamma(C, H, T)e^{\langle \la, H \rangle}q(H)\,dH, \quad \la \in V_{\C}.
\end{equation}
Lemma \ref{lem:gamma1}, and references invoked in the beginning of the paragraph, 
ensure that $\scrF(\Gamma(C), T, q,\la)$ is defined by an absolutely convergent integral. 
Moreover, an analysis similar to the one effectuated in \cite{arthur2}, Lemma 2.2 or \cite{leMoi}, Lemme 4.3, 
shows that for a fixed $\la \in V_{\C}$, the function
\[
T \in V \to \scrF(\Gamma(C), T,q,\la)
\] 
is a polynomial exponential. Additionally, 
for $\la \in V_{\C}^{C-reg}$
the purely polynomial part of $\scrF(\Gamma(C), T,q,\la)$ is constant and given by
$\scrF(C, q,\la)$.

\section{Absolute theory}\label{sec:redTh}
\subsection{Generalities}\label{ssec:generalities}

Let $F$ be a number field. 
Let $\A$ be the ring of adeles of $F$. Let $|\cdot|_{\A}$ be the 
norm on the ideles $\A^{\times}$ of $F$ normalized in the standard way. 
We write $\A = F_{\infty} \times \A^{\infty}$ where $F_{\infty} = F \otimes_{\Q}\R $  
and $\A^{\infty}$ is the ring of finite adeles. 
 Let $G$ be a connected reductive algebraic group defined over $F$.
All subgroups of $G$ are assumed to be closed and defined over $F$.

Let $\Gm$ be the multiplicative group 
and $\Ga$ - the additive group, both seen as linear algebraic groups defined over $F$. 
For a subgroup $H \sbs G$, let $X_{*}(H) = \Hom_{F}(\Gm, H)$ and $X^{*}(H) = \Hom_{F}(H,\Gm)$.

Fix $A_0$ a maximal $F$-split torus of $G$.  
Parabolic subgroups of $G$ containing $A_0$ are called semi-standard 
and their set is denoted $\calF(A_{0})$. In general, for a subgroup $H$ of $G$
we denote $\calF^{G}(H) = \calF^{G}(H) $ the set of semi-standard parabolic subgroups of $G$ 
containing $H$. 

If $P \in \calF(A_{0})$ is a parabolic subgroup we denote $N = N_{P}$ 
its unipotent radical 
and $M = M_{P}$ its unique Levi subgroup containing $A_{0}$. 
We denote also $A_{P}$ be the split center of $M = M_{P}$. 
It can also be defined as the subtorus of $A_{0}$ centralizing $M$. 
For a subtorus $A_1 \sbs A_{0}$ we set $\calP(A_{1})$ to be the set of semi-standard parabolic subgroups $P$
such that $A_{P} = A_{1}$.
We reserve the letters $P, Q, R$ and sometimes $S$ to denote (semi-standard) parabolic subgroups.

\subsection{Geometry of chambers}\label{ssec:chambs}

Let $\ago_{0} = X_{*}(A_0) \otimes_{\Z} \R$.
For $P \in \calF(A_0)$ we also put $\ago_{P} = X_{*}(A_P) \otimes_{\Z} \R$. 
The space $\ago_{P}$ only depends on the Levi part $M$ of $P$ 
and will sometimes be denoted by $\ago_{M}$. We have naturally $\ago_{P} \sbs \ago_{0}$.
Fix a $W$-invariant scalar product $\bilif$ on $\ago_{0}$. 
Thanks to it, the groups $X^{*}(A_{P})$ and $X^{*}(M_P)$ are naturally realized as lattices in $\ago_{0}$, 
contained in $X_{*}(A_0) \otimes_{\Z} \Q = X^{*}(A_0) \otimes_{\Z} \Q$, 
and $\ago_{P}^{*} := X_{*}(A_{P}) \otimes_{\Z} \R$ gets identified with $\ago_{P}$. 

The scalar product induces a topology on $\ago_{0}$. We write $\bar A$ 
for a closure of a subset $A \sbs \ago_{0}$ in this topology.

Let $P,Q \in \calF(A_{0})$ be such that $P \sbs Q$. 
Let $\ago_{P}^{Q}$ be the orthogonal complement of $\ago_{Q}$ in $\ago_{P}$ 
and $\Delta_{P} \sbs \ago_{P}$ be the set of simple roots 
for the action of $A_{P}$ on $N_{P}$. 
Let $\Delta_{P}^{Q} \sbs \ago_{P}$ be the subset of $\Delta_{P}$ that vanishes 
on $\ago_{Q}$. It is naturally contained in $\ago_{P}^{Q}$, where it forms a basis.
We define hence $\hDelta_{P}^{Q} \sbs \ago_{P}^{Q}$ as the dual basis to $\Delta_{P}^{Q}$.
If $Q =G$, we write simply $\Delta_{P} = \Delta_{P}^{G}$. 

Define also
\[
\ago_{P}^{+} = \{ H \in \ago_{P} \ | \ \langle H, \al \rangle > 0 \ \forall \al \in \Delta_{P}\}.
\]
If $Q =G$, we write simply $\Delta_{P} = \Delta_{P}^{G}$, $\ago_{P}^{ +} = \ago_{P}^{G, +}$ etc.
We can and will use the language introduced in Section \ref{sec:cones} to work with the above chambers. 
First of all, by definition, we have $\rint \bragop = \ago_{P}^{+}$. 
The cones $\bragoq$ run over faces of the cone $\bragop$ for $Q \sps P$. 
Moreover, we have in this case
\begin{equation}\label{eq:angRoot}
A(\bragoq, \bragop) = \{ H \in \ago_{P} \ | \ \langle H, \al \rangle \ge  0 \ \forall \al \in \Delta_{P}^{Q}\}.
\end{equation}

The following lemma lists standard facts about root chambers whose proofs can be found in Paragraph 1.2 of \cite{labWal}.

\blem\label{lem:basicRootProps} Let $P,Q \in \calF(A_{0})$ be such that $P \sbs Q$.  The following properties hold. 
\begin{enumerate}
\item For all $\al, \al' \in \Delta_{P}^{Q}$ we have $\langle \al, \al' \rangle \le 0$.
\item For all $\varpi, \varpi' \in \hDelta_{P}^{Q}$ we have $\langle \varpi, \varpi' \rangle \ge 0$.
\item $A(\bragoq, \bragop)\cap \ago_{P}^{Q} \sbs A(\bragoq, \bragop)^{\vee}$.
\item The projection of $\ago_{P}^{+}$ onto $\ago_{Q}$ 
is contained in $\ago_{Q}^{+}$.
\end{enumerate}
\elem

The next result is essentially Lemme 1.7.1 of \cite{labWal}.

\blem\label{lem:varpivarpi} Let $P_0 \in \calP(A_0)$ and $P \in \calF(P_{0})$. Then, 
there exist an $m \ge 0$ 
such that for all $H \in \ago_{0}$ and all $C_1,C_2 \in \R$
such that 
\[
\langle \varpi, H \rangle \ge  C_1, \ \forall \ \varpi \in \hDelta_{P}, \quad 
\langle \overline  \varpi, H \rangle  \ge C_2, \ \forall \ \overline \varpi \in \hDelta_{0}^{P} 
\]
we have 
\[
\langle \varpi, H \rangle \ge C_2 - m |C_1|, \quad \forall \ \varpi \in \hDelta_{0} \smin \hDelta_{P}.
\]
\elem
\bdem 
Let $\varpi \in \hDelta_{0} \smin \hDelta_{P}$. We can write it uniquely as $\varpi = \overline \varpi + \varpi'$
where $ \overline \varpi \in \hDelta_{0}^{P}$ and $\varpi' \in \ago_{P}$. 
Using part 4. of Lemma \ref{lem:basicRootProps} we have $\varpi' \in \bragop$ 
and the result follows. 
\edem

\subsection{The $\tlGam$ function}\label{ssec:tlgam}

Let $P,R \in \calF(A_0)$ such that $R \sbs P$. We define the following function 
of variables $H, T \in \ago_{0}$
\[
\tlGam_{R}^{P}(H,T) := \Gamma(A(\bragop, \bragor), H, T).
\]
This is essentially the function $\Gamma'_{R \cap M}$ (associated to the group $M$ 
and its parabolic subgroup $R \cap M$) studied by Arthur \cite{arthur2}. 
We will rather use \cite{labWal}, Paragraph 1.8, as reference where this function 
is noted $\Gamma^{M}_{R \cap M}$.
We reserve however the notation $\Gamma_{R}^{P}$ for further purposes. 

For $H \in \ago_{0}$ let $H = H_P + H^{P}$ 
be its decomposition with respect to the orthogonal sum 
decomposition $\ago_{0} = \ago_{P} \oplus \ago^{P}$.

\blem\label{lem:oldGamma} Let $T \in \ago_{0}$. 
The function 
\[
H \in \ago_{0} \to \tlGam_{R}^{P}(H,T)
\]
is compactly supported. Moreover, if $T \in \ago_{0}^{+}$, it is a characteristic function of 
$H \in \ago_{0}$ such that
\begin{enumerate}
\item $H_{P} = T_P$;
\item $H \in \ago_{R}$;
\item $\langle \al, H \rangle > 0$ for all $\al \in \Delta_{R}^{P}$;
\item $\langle \varpi, H - T\rangle \le 0$ for all $\varpi \in \hDelta_{R}^{P}$.
\end{enumerate}
\elem
\bdem
The first two points follow from Lemma \ref{lem:gamma0}. 
Suppose then that $H_{P} = T_P$ and $H \in \ago_{R}$.
We have 
\[
\tlGam_{R}^{P}(H,T) = \Gamma^{M}_{R \cap M}(H^{P}, T^{P})
\]
where $\Gamma^{M}_{R \cap M}$ is studied in \cite{labWal}, Paragraph 1.8. 
The desired result is proven therein, Lemme 1.8.3 (with "r\'egulier" defined at the end of Paragraph 1.2 
as belonging to $\ago_{0}^{+}$ ). 
\edem 

\brem As is clear from the lemma above, we slightly deviate from the standard 
usage of the function $\tlGam_{R}^{P}(\cdot, T)$. Classically, it is studied as a function on 
$\ago_{P}^{R}$. We insist however on considering it as a function on the whole space $\ago_{0}$. This point of view is in line 
with our approach from Section \ref{sec:cones}. 
\erem

The next lemma is proved in \cite{labWal}, Lemme 1.2.7

\blem\label{lem:labWal127} Let $P \sbs Q$. Suppose $\langle \al, H \rangle >0 $
for all $\al \in \Delta_{P}^{Q}$ 
and $\langle \varpi, H\rangle  \le 0$ for all $\varpi \in \hDelta_{0}^{P}$. 
Then $\langle \al, H \rangle >0 $
for all $\al \in \Delta_{0}^{Q} \smin \Delta_{0}^{P}$.
\elem

\subsection{Parabolic subgroups}\label{ssec:springCoch} 

We invoke the fact, \cite{springer} Paragraph 13.4, that parabolic subgroups 
of $G$ are obtained from cocharacters. 
Recall first that for a morphism $\phi : \Gm \to X$ of algebraic varieties defined over $F$,
we say $\lim_{a \to 0} \phi(a)$ exists, if $\phi$ extends to an algebraic 
morphism $\tlphi : \Ga \to X$. We put then $\lim_{a \to 0} \phi(a) = \tlphi(0)$.

We have that $\la \in X_{*}(A_{0})$ 
gives rise to the semi-standard parabolic subgroup
\[
P(\la) := \{x \in G \ | \ \lim_{a \to 0}\la(a)x\la(a)^{-1} \text{ exists} \}.
\]
Then, $P = P(\la)$ if and only if $\la \in X_{*}(A_{0}) \cap \ago_{P}^{+}$.

Suppose then that $\la \in X_{*}(A_{0}) \cap \ago_{P}^{+}$. 
Then
\[
N = \{x \in P \ | \ \lim_{a \to 0}\la(a)x\la(a)^{-1} = e \}
\]
where $e$ is the identity element of $G$. Moreover, the Levi part $M$ of $P$ is the centralizer of the image of $\la$ 
in $G$.

\subsection{Reduction theory}\label{ssec:redTh}

Let $W$ be the Weyl group associated to the couple $(G,A_{0})$.
We choose representatives of elements of $ W $ to lie in $ G(F) $. 
We take a maximal compact $ K $ of $ G(\A) $, adapted to $ M_{0} $. 
Let $P \in \calF(A_{0})$. Define $A_{P}^{\infty}$ 
as the connected component of the $\R$-points of $\Res_{F/\Q}A_{P}$ - the restriction 
of scalars of $A_{P}$ from $F$ to $\Q$. 
Define the Harish-Chandra function $H_{P} : G(\A) \to \ago_{P}$ 
so that the following relation is satisfied
\[
\langle \xi, H_{P}(x) \rangle = \log |\xi (m)|_{\A}
\]
for all $\xi \in X^{*}(M)$, where we decompose $x = nmk$ using the Iwasawa decomposition $G(\A) = N(\A)M(\A)K$.

In particular, $H_{G}$ is a homomorphism and we denote $G(\A)^{1}$ its kernel. 
We have then the direct product decomposition of commuting subgroups $G(\A) = G(\A)^{1}A_{G}^{\infty}$.

We define a height $|\cdot| : G(\A) \to \R_{+}^{\times}$ as in \cite{labWal}, Paragraph 3.2. 
It is a continuous map such that there exist positive constants $c_1$, $c_2$
such that for all $x,y \in G(\A)$ we have 
\begin{itemize}
\item[-] $|x| \ge c_1$;
\item[-] $|xy| \le c_2 |x||y|$;
\item[-] $|x| = |x^{-1}|$.
\end{itemize}

%

Fix $P_{0} = M_{0}N_{0} \in \calF(A_{0})$, a minimal Levi subgroup of $ G $. 
Write $\Delta_{0}$ for $\Delta_{P_{0}}$, $H_0$ for $H_{P_{0}}$ etc.
For a compact $ \omega \sbs P_{0}(\A)$, $ T \in \ago_{0}^{G} $ 
and $P \in \calF(P_0)$
let
\[ 
\Sgl(P,T, \omega)= \{ cak| c \in \omega, \ k \in K,  \ a \in A_{0}^{\infty}, \ \langle \al, H_{0}(a) - T \rangle \ge 0,
\, \forall \ 
\al \in \Delta_{0}^{P}\}.
\]
We can then choose, and we do, $ \omega $ and $ T = \tzers $ so that 
$ G(\A) = P(F)\Sgl(P, \tzers, \omega) $ for all $P \in \calF(P_0)$.
For this choice of $\tzers$ and $\omega$ we simply set 
$\Sgl(P) = \Sgl^{G}(P)  = \Sgl(P, \tzers, \omega)$.

Let us recall some classical results from reduction theory. 
We denote $\|\cdot\|$ the norm on the Euclidean space $\ago_{0}$. 

\blem\label{lem:redTh} There exist constants $c_4$, $c_5$ and $c_6$ 
such that for all $x \in G(\A)$ the following holds
\begin{enumerate}
\item 
\[
\|H_0(x)\| \le c_4(1 + \log |x|).
\]
\item There exists a $\delta \in P_0(F)$ such that 
\[
\log |\delta x| \le c_5(1 + \|H_0(\delta x)\|).
\]
\item If $x \in \Sgl(G)$ then for all $\gamma \in G(F)$ and all $\varpi \in \hDelta_{0}$ we have
\[
\varpi(H_{0}(\gamma x)) \le \varpi(H_{0}(x)) + c_6.
\]
\end{enumerate}
\elem
\bdem
All three assertions are proven in \cite{labWal}, Lemmas 3.2.2, 3.5.3 and 3.5.4 respectively. 
\edem

The proof of the following lemma follows closely the proof of Lemma 3.6.1 of \cite{labWal}.

\blem\label{lem:lang212}  There exists a $\tregs \in \ago_{0}$ 
with the following property. For all $T \in \tregs + \ago_{0}^{+}$, all $x \in \Sgl(G)$ such that 
\[
\al(H_{0}(x)) > \al(T) \quad \forall \ \al \in \Delta_{0} \smin \Delta_{0}^{P}
\]
and all $\delta \in G(F)$ we have that $\delta x \in \Sgl(G)$ implies that $\delta \in P(F)$.
\elem

\bdem 

Suppose that $\al(H_{0}(x)) > \al(T)$ for all $\al \in \Delta_{0} \smin \Delta_{0}^{P}$ and some $T \in \ago_{0}^{+}$ 
and that $\delta x \in \Sgl(G)$ for some $\delta \in G(F)$. 
We will show that there exists a constant $C$, that depends only on $G$, 
such that if $\al(T) > C$ for all $\al \in \Delta_{0}$, then we must have $\delta \in P(F)$. 
This clearly implies the desired statement. 

Using the Bruhat decomposition we can replace $\delta$ with $w_{s}$ - 
an element in $G(F)$ representing an element $s$ of the Weyl group $W$. 
Our goal is thus to show that $s$ belongs to the Weyl group of $M$.
In the following, by an "absolute constant", we mean a constant that depends only on $G$, and not on $x$ 
or $T$.

Let $\Sigma^{+}(P_{0})$ be the set of positive roots with respect to $A_0$ and $P_0$. 
Define
\[
\calR(s) = \{ \beta \in \Sigma^{+}(P_{0}) \ | \ s\beta \nin \Sigma^{+}(P_{0})\}.
\]
Let $\rho_{0}$ be the half of the sum of elements of $\Sigma^{+}(P_{0})$ and
\[
\la = \rho_{0} - s^{-1}\rho_{0}.
\]
Note that 
\begin{equation}\label{eq:redTh0}
\la = \sum_{\beta \in \calR(s)}\beta.
\end{equation}
We claim that there exists an absolute  constant $C'$ such that
\begin{equation}\label{eq:redTh}
\langle \la, H_{0}(x) \rangle \le C'.
\end{equation}

Let's prove \eqref{eq:redTh}. 
By Lemma 3.3.2 of \cite{labWal}, we have
\begin{equation}\label{eq:redTh1}
H_0(x) = s^{-1}H_{0}(w_{s}x) - \sum_{\beta \in \calR(s)}c_{\beta}\beta
\end{equation}
where the coefficients $c_{\beta}$ satisfy $c_{\beta} \ge -c$ for an absolute constant $c$. 

For all $\beta \in \calR(s)$ we have 
\begin{equation}\label{eq:redTh2}
\langle \la, \beta \rangle = \langle \rho_{0}, \beta + (-s\beta) \rangle > 0.
\end{equation}
On the other hand, since 
\[
s\la = -\sum_{\gamma \in \calR(s^{-1})}\gamma
\]
and 
\[
\langle \beta, H_{0}(w_{s} x) - \tzers \rangle \ge 0,\quad \forall \ \beta \in \Sigma^{+}(P_{0})
\]
by assumption $w_{s}x \in \Sgl(G)$, 
we have
$
\langle s\la, H_{0}(w_{s} x) - \tzers \rangle \le 0
$
which means 
\begin{equation}\label{eq:redTh3}
\langle \la, s^{-1}H_{0}(w_{s} x) \rangle \le c_1
\end{equation}
for come absolute constant $c_1$.

Combining \eqref{eq:redTh1}, \eqref{eq:redTh2} and \eqref{eq:redTh3} we obtain \eqref{eq:redTh}. 

Using the decomposition of $\la$ into a sum of roots \eqref{eq:redTh0} and the assumption 
on $H_{0}(x)$ we see that \eqref{eq:redTh} cannot be maintained for $\al(T)$ large 
unless all elements of $\calR(s)$ are contained in $\Sigma^{+}(P_0 \cap M)$. 
This means in turn that $s$ lies in the Weyl group of $(M, A_0)$.
\edem

The next Lemma is a generalization of Lemma 5.1 of \cite{arthur3}. 
See also Lemme 3.7.1 of \cite{labWal}.

\blem\label{lem:htauFinite} Let $P \in \calF(P_0)$, $\la \in \ago_{P}^{+}$ and $T \in \ago_{0}$. Then, there exists a constant $C > 0$ 
and $N \in \N$ such that for all $x \in \Sgl(G)$
the number of elements in the set
\[
\{   \delta \in P(F) \bsl G(F) \ | \ \langle H_0(\delta x) - T, \la \rangle > 0 \}
\]
is finite and bounded by
\[
C(|x|e^{\|T\|})^{N}.
\]
\elem
\bdem
Let's write
\[
\la = \sum_{\varpi \in \hDelta_{P}}a_{\varpi}\varpi + \la_{G}
\]
where $\la_{G} \in \ago_{G}$ and $a_{\varpi} \in \R_{+}^{\times}$ by assumption. 
Let  $y \in G(\A)^{1}$ and $z \in A_{G}^{\infty}$ be such that $x = yz$.

Let $\delta \in P(F) \bsl G(F)$ be such that $\langle H_0(\delta x) - T, \la \rangle > 0$. 
By assumption, we have
\begin{equation}\label{eq:theAss}
\langle \la, H_{0}(\delta y) \rangle = 
\sum_{\varpi \in \hDelta_{P}}a_{\varpi}\langle \varpi, H_{0}(\delta y) \rangle
> \langle \la, T - H_{G}(z) \rangle. 
\end{equation}
Let us put $C_1 := \langle \la, T - H_{G}(z) \rangle$.

On the other hand, part 3 of Lemma \ref{lem:redTh} says that there exists a constant $c_6$ such that
\[
\langle \varpi, H_{0}(\delta y) \rangle \le \langle \varpi, H_{0}(y) \rangle + c_6, \quad \forall \ \varpi \in \hDelta_{0}.
\]
This implies that there is a constant $c$ such that 
\begin{equation}\label{eq:theHwn}
\langle \varpi, H_{0}(\delta y) \rangle  \le c(\|H_{0}(y)\| + 1), \quad \forall \ \varpi \in \hDelta_{0}.
\end{equation}
Let us put $C_2 = c(\|H_{0}(y)\| + 1)$.

Let $a = \max_{\varpi \in \hDelta_{P}} a_{\varpi}$ and $b = \min_{\varpi \in \hDelta_{P}} a_{\varpi}$.
Put
\[
C_3 = \dfrac{1}{b}(-|C_1| - a(|\hDelta_{P}|-1)C_2)
\]
where $|\hDelta_{P}|$ denotes the number of elements in $\hDelta_{P}$.

We claim that for all $\varpi \in \hDelta_{P}$ we have 
\begin{equation}\label{eq:theClm}
\langle \varpi, H_{0}(\delta y) \rangle > C_3.
\end{equation}
Indeed, let us prove by reduction to contradiction. Suppose 
there is a $\varpi' \in \hDelta_{P}$ such that 
\[
\langle \varpi', H_{0}(\delta y) \rangle \le C_3.
\]
We add this inequality multiplied by $a_{\varpi'}$ 
to the $|\hDelta_{P}|-1$ inequalities \eqref{eq:theHwn}  for all $\varpi \in \hDelta_{P} \smin \{\varpi'\}$
each multiplied by the corresponding factor $a_{\varpi}$. 
We obtain hence
\[ 
\langle \la, H_{0}(\delta y) \rangle \le 
a_{\varpi'}C_3 + 
\sum_{\varpi \in  \hDelta_{P} \smin \{\varpi'\}}
a_{\varpi} C_2 \le bC_3 + a(|\hDelta_{P}|-1)|C_2 = -|C_1|
\]
which contradicts \eqref{eq:theAss}.

We can assume that $\delta x \in \Sgl(P)$ 
and that the conditions of Lemma \ref{lem:redTh}, part 2, are satisfied. 
Since $\delta x \in \Sgl(P)$ we also have $\delta y \in \Sgl(P)$. 
That means in particular
\[
\langle \al, H_{0}(\delta y) - \tzers \rangle > 0 \quad \forall \ \al \in \Delta_{0}^{P}.
\]
Using part 3 of Lemma \ref{lem:basicRootProps} we obtain thus
\[
\langle \overline \varpi, H_{0}(\delta y) - \tzers \rangle > 0  \quad \forall \ \overline \varpi \in \hDelta_{0}^{P}.
\]

The lemma now easily follows from Lemma \ref{lem:varpivarpi} and parts 1 and 2 of Lemma \ref{lem:redTh} 
(see Lemme 3.7.1 of \cite{labWal}).
\edem

We fix a $\tregs \in \ago_{0}$ satisfying the conditions of Lemma \ref{lem:lang212}. 

\subsection{Estimates}\label{ssec:altSum}

For a linear algebraic group $H$, let $\Ugo(H)$ to be the envelopping algebra 
of the complexification of $\Lie(H)(F_{\infty})$.
Its elements act on the right on $C^{\infty}(H(F_{\infty}))$. 
We note this action $\rho(X) f$ for $f \in C^{\infty}(H(F_{\infty}))$
and $X \in \Ugo(H) $.

 Let $P \in \calF(A_{0})$. For $\al \in \Delta_{P}$ 
Let $N_{P, \al}$ be the unipotent radical of the maximal parabolic subgroup $P_{\al}$ 
such that $\Delta_{P} \smin \Delta_{P}^{P_{\al}} = \{\al\}$.

For a smooth function  $\phi$ on $N_{P}(F) \bsl N_{P}(\A)$ define
\[
\phi_{P,Q}= \sum_{P \sbs R \sbs Q}
(-1)^{\dim \ago_{R}/\ago_{Q}}\phi_{R}
\]
where
\[
\phi_{R}(x) = \int_{[N_{R}]}\phi(nx)\,dn.
\]

We have the following result proved in \cite{labWal}, Lemme 4.3.1.

\blem\label{lem:431} 
Fix an open compact subgroup $\calO$ 
of $N(\A_{\infty})$. 
For all $r \ge 0$ there exists a finite number of operators $X_{k} \in \Ugo(N_{P})$
of the form
\[
X_{k} = \prod_{\al \in \Delta_{P}^{Q}} \displaystyle \left (   \sum_{j=1}^{n_{P, \al}^{Q}} Y_{k,\al,j}^{r}  \right ),\quad 
Y_{k,\al, j} \in \Lie(N_{P,\al} \cap M_{Q})(F_{\infty})
\]
where $n_{P, \al}^{Q} = \dim_{F}(N_{P,\al} \cap M_{Q})$, such that
\[
\sup_{n \in [N_{P}]} |\phi_{P,Q}(n)| \le \sum_{k}
\sup_{n \in [N_{P}]} |\rho(X_{k}) \phi(n)|
\]
for all smooth function $\phi$ on $N_{P}(F) \bsl N_{P}(\A) / \calO$.
\elem

We also note the following elementary identity
\begin{equation}\label{eq:constAlt}
\phi_{Q} = \sum_{P \sbs R \sbs Q}\phi_{R}\sum_{R \sbs S \sbs Q}
(-1)^{\dim \ago_{R}/\ago_{S}} = 
\sum_{P \sbs S \sbs Q} \sum_{P \sbs R \sbs S}
(-1)^{\dim \ago_{R}/\ago_{S}} \phi_{R} = 
\sum_{P \sbs S \sbs Q} \phi_{P,S}
\end{equation}

\section{Relative theory}\label{sec:relRedTh}
\subsection{A subgroup of $G$}\label{ssec:sbgp}

Let $G'$ be a connected reductive subgroup of $G$.
Fix $A_{0}'$  a maximal $F$-split torus of $G'$.
Fix also $P_{0}'= M_{0}' N_{0}' \in \calP^{G'}(A_{0}')$ a minimal parabolic subgroup of $G'$. 
Fix $K'$ a good maximal compact of $G'(\A)$ adapted to $M'_{0}$, etc. 
We use the results of Section \ref{sec:redTh} applied to the group $G'$. 
All data attached to $G'$ will be denoted with a prime.

Let $M_{1}$ be the centralizer of $A_{0}'$ in $G$. It is a Levi-subgroup 
of $G$. Let $A_{1}$ be its split center. 
We can, and do, assume that $A_{0} \sps A_{1}$
We have then $\ago_{0'} \sbs \ago_{1} \sbs \ago_{0}$. 

For $P' \in \calF^{G'}(A_{0}')$ set
\[
\calP^{G}(P') = \{ P \in \calF^{G}(A_{1}) \ | \ \ago_{P'}^{+} \cap \ago_{P}^{+}  \neq \varnothing\}.
\]
\brop\label{prop:PcapH} 
For all $P \in \calP^{G}(P')$ we have $P \cap G' = P'$. In particular, we 
have 
\[
 \ago_{P'}^{+} \cap \ago_{P}^{+}  = \ago_{P'} \cap \ago_{P}^{+} =  \ago_{0'} \cap \ago_{P}^{+}.
 \]
Moreover, $M \cap G' = M'$ and $N \cap G' = N'$.
\erop
\bdem

The spaces $\ago_{P}$ and $\ago_{P'}$ come with obvious, compatible rational structures. 
The cones $\ago_{P}^{+}$ and $\ago_{P'}^{+}$ 
are open in $\ago_{P}$ and $\ago_{P'}$. By assumption $\ago_{P}^{+} \cap \ago_{P'}^{+}$
is open and non-empty in $\ago_{P} \cap \ago_{P'}$. Invoking Proposition \ref{prop:brunGub} 
we see that $\ago_{P}^{+} \cap \ago_{P'}^{+}$ is given by  intersection 
of rational half-spaces 
and therefore must have a rational, and consequently also a lattice point. 
An element of 
\[
X_{*}(A_{0}') \cap \ago_{P'}^{+} \cap \ago_{P}^{+} \neq \varnothing
\]
demonstrates the desired statements by the results recalled in \ref{ssec:springCoch}. 
\edem

We set as well
\[
\calF^{G}(P') = \bigcup_{Q' \in \calF^{G'}(P')}\calP^{G}(Q').
\]

\brem Note that if $A_{1} = A_{0'}$ (for example if the split ranks of $G$ and $G'$ are equal)
then for any $P \in \calF^{G}(A_{1})$ we have that $P \cap G'$ is a semi-standard 
parabolic subgroup of $G'$ and $P \in \calF^{G}(P_{0}')$ if and only if $P_{0} \sbs P \cap G' $.
\erem

\subsection{Relative chambers}

We put for $P \in \calF^{G}(P_{0}')$
\[
\zgo_{P} := \ago_{P} \cap \ago_{0'} , \quad \zgo_{P}^{+} := \ago_{P}^{+} \cap \ago_{0'}.
\]
The cone $\zgo_{P}^{+}$ is an open cone in $\zgo_{P}$. 
We have disjoint union decompositions for $P' \in \calF^{G'}(P_{0})$
\begin{equation}\label{eq:bigFan}
\ago_{P'}^{+} = \bigsqcup_{P \in \calP^{G}(P')}\zgo_{P}^{+}, \quad 
\overline{\ago_{P'}^{+}} = \bigsqcup_{P \in \calF^{G}(P')}\zgo_{P}^{+}.
\end{equation}
Let $P \in \calF^{G}(P_{0}')$.
The map
\[
Q \in \calF^{G}(P) \cap \calF^{G}(P_{0}') \mapsto \brzgoq
\]
is a bijection between $ \calF^{G}(P) \cap \calF^{G}(P_{0}')$ and $\calF(\brzgop)$, i.e. the set 
of faces of the cone $\brzgop$. 
In particular, if $P \neq G$, the cone $\brzgop$ 
is not a linear space. 
 
We are going to use the results of Section \ref{sec:cones}
with $V = \ago_{0'}$. The Euclidean, Weyl group invariant structure on $\ago_{0'}$ 
is that induced from the scalar product $\bilif$ on $\ago_{0}$. 
Since such a structure can be defined using an embedding into some $GL_n$, clearly they can be chosen 
in a compatible way.
For $P,Q \in \calF^{G}(P_{0}')$ such that $P \sbs Q$
define 
\begin{itemize}[-]
\item $\tau_{P}^{Q}$ - the characteristic function of $\rint A(\brzgoq, \brzgop)$.
\item $\htau_{P}^{Q}$ - the characteristic function of $\rint A(\brzgoq, \brzgop)^{\vee}$.
\item $\zgo_{P}^{Q}$ - the orthogonal complement of $\zgo_{Q}$ in $\zgo_{P}$.
\item $\varepsilon_{P}^{Q} = (-1)^{\dim \zgo_{P}^{Q}}$.
\item $X_{P}$, $X^{P}$, $X_{P}^{Q}$ - projections of an element $X \in \ago_{0'}$ 
onto $\zgo_{P}$, $\zgo^{P}$ and $\zgo_{P}^{Q}$ respectively.
\end{itemize}

\subsection{The $\sigma$ and $\Gamma$ functions} 

Let $P,Q \in \calF^{G}(P_{0}')$ such that $P\sbs Q$. Set
\[
\Gamma_{P}^{Q}(H,X) := \Gamma(A(\brzgoq, \brzgop), H,X), \quad H,X \in \ago_{0'}
\]
where the right hand side function is defined in Section \ref{ssec:gamma}. Set also
\[
\sigma_{P}^{Q}(H) := \sigma(\brzgoq, \brzgop)(H), \quad H \in \ago_{0'}
\]
where the right hand side function is defined in Section \ref{ssec:sigma}.

\blem\label{lem:conesProps} 
\begin{enumerate}[1.]
\item For all $X \in \ago_{0'}$, the function $\Gamma_{P}^{Q}(\cdot,X)$ 
is compactly supported.
\item 
\[
\tau_{P}^{Q}(H) = 
\ssum{
R \in \calF^{G}(P_{0}')\\
P \sbs R \sbs Q 
} \Gamma_{P}^{R}(H^{R},X^{R})\tau_{R}^{Q}(H_{R}-X_{R}) \quad \forall \, H,X \in \ago_{0'}.
\]
\item 
\[
\htau_{P}^{Q}(H-X) = 
\ssum{
R \in \calF^{G}(P_{0}') \\
P \sbs R \sbs Q 
} \varepsilon_{R}^{Q} \htau_{P}^{R}(H^{R}) \Gamma_{R}^{Q}(H_{R},X_{R}), \quad \forall \, H,X \in \ago_{0'}.
\]
\item 
\[
\htau_{Q}\tau_{P}^{Q} = 
\ssum{
R \in \calF^{G}(P_{0}') \\
R \sps Q 
} \sigma_{P}^{R}.
\]
\item 
We have for $H \in \ago_{0}'$
\[
\sigma_{P}^{Q}(H) \neq 0 \Longrightarrow \tau_{P}^{Q}(H) \neq 0.
\]
Moreover, there exist a constant $k > 0$ such that 
for all $H \in \ago_{0'}$ such that $\sigma_{P}^{Q}(H) \neq 0$ 
we have
\[
\|H_{Q}\| \le k\|H^{Q}\|.
\]
\end{enumerate}
\elem
\bdem
The first point is Lemma \ref{lem:gamma1}, the second Lemma \ref{lem:gamma2}, the third 
Corollary \ref{cor:gamma2}, the fourth Lemma \ref{lem:htauTauSigSum} and the fifth Proposition \ref{prop:sigma}. 
\edem

\subsection{Relative decomposition of $M_{P'}(\A)$}

Let $P, Q \in \calF^{G}(P_{0}')$ be such that $P \sbs Q$.
Let $P' \in \calF^{G'}(P_{0})$. For $x \in G'(\A)$ we will write
\[
H_{P'}(x)_{P} = (H_{P'}(x))_{P}, \quad H_{P'}(x)^{Q} = (H_{P'}(x))^{Q}, \quad  H_{P'}(x)^{Q}_{P} = (H_{P'}(x))^{Q}_{P}.
\]

Suppose $P'= P \cap G$. Let
\[
M_{P'}(\A)^{P,1} := \{x \in M_{P'}(\A) \ | \ H_{P'}(x)_{P} = 0\}.
\]
Set $Z_{P}^{\infty} = A_{P}^{\infty} \cap A_{P'}^{\infty}$. 
The restriction of $H_{P'}$ to $Z_{P}^{\infty}$ is a group isomorphism with its image $\zgo_{P}$.
We obtain hence the direct 
product decomposition of commuting groups
\[
M_{P'}(\A) = Z_{P}^{\infty}M_{P'}(\A)^{P,1}.
\]

\subsection{The $F^{P}$ function}
First let us make the following notation. 
\[
\tzer := T_{G'}, \quad \treg := T_{G',reg}.
\]
These are elements of $\ago_{0'}$. The element $T_{G'}$ was defined in Paragraph \ref{ssec:redTh} 
and $T_{G', reg}$ is such that Lemma \ref{lem:lang212} holds for $G'$. 

Let $P \in  \calF^{G}(P_{0}')$ and $P' = P \cap G'$. 
For $T \in \ago_{0'}$ let 
$F^{P}(\cdot , T)$ be the characteristic function of $x \in G'(\A)$
such that there exists a $\delta \in P'(F)$
such that $\delta x \in \Sgl(P')$ 
and
\begin{equation}\label{eq:FPrel}
\tlGam_{0'}^{P'}(H_{0'}(\delta x)^{P} - \tzer^{P}, T^{P} - \tzer^{P}) = 1
\end{equation}
where $\tlGam_{0'}^{P'}$ is defined in Paragraph \ref{ssec:tlgam}.

Some comments are in place. Our function $F^{P}$ is inspired by Arthur's, defined in \cite{arthur3}, Section 6. 
It is often referred to as the characteristic function of the truncated fundamental domain. 
Let us compare briefly Arthur's construction and ours. First of all, the truncation takes place at the level of the group $G'$. 
Therefore, $F^{P}$ should be compared with Arthur's $F^{P'}$. Our $F^{P}$ isn't really related 
to the Arthur's $F^{P}$ defined from the point of view of the group $G$. 

The difference between $F^{P}$ and $F^{P'}$ is in fact quite subtle. Arthur defines $F^{P'}(\cdot, T)$ requiring $\delta x$ to satisfy
\begin{equation}\label{eq:FPrel2}
\tlGam_{0'}^{P'}(H_{0'}(\delta x)^{P'} - \tzer^{P'}, T^{P'} - \tzer^{P'}) = 1.
\end{equation}
Comparing this to \eqref{eq:FPrel} we see that the difference lies in the way we project the vectors. 
We have that $\ago_{0'}^{P'}$ is contained in the orthogonal complement of $\zgo_{P}$ in $\ago_{0'}$. 
Therefore, we simply have that \eqref{eq:FPrel} holds if and only if \eqref{eq:FPrel2} holds 
and 
\[
H_{P'}(\delta x)^{P} = T_{P'}^{P}.
\]
For some $P$ the above condition is trivial, for other it is not. In particular, 
if $A_{G'} \neq A_{G}$, the function $F^{G}(\cdot ,T)$, defined 
on $G'(F) \bsl G'(\A)^{G,1}$ 
is supported on a set of measure $0$. 
Figure \ref{figure:4} describes this eventuality. 
See also comments after Proposition \ref{prop:decOf1Rel} and its proof for further comments on this.

\subsection{Relative decomposition of $1$}
Let $Q \in \calF^{G'}(P_0')$ and note $Q'= Q \cap G'$. The following result is 
the generalization of the classical result, Lemma 6.4 of \cite{arthur3} to the relative case.
Our proof combines the proof of loc. cit. with the methods used in \cite{labWal}, Proposition 3.6.3.

\brop\label{prop:decOf1Rel}
 For all $x \in G'(\A)$ and all $T \in \treg + \ago_{0'}^{+}$ we have
\[
 \ssum{
 P \in \calF^{G}(P_{0}') \\
 P \sbs Q 
}
\sum_{\delta \in (P \cap G')(F) \bsl Q'(F)} 
F^{P}(\delta x, T)\tau_{P}^{Q}(H_{0'}(\delta x)_{P} - T_{P}) = 1.
\]
\erop 

Let us give some context before the proof. Up to reduction theory, 
such statement reduces to combinatorics of cones. 
In the classical case, for the group $G'$, 
the associated combinatorial result is Lemma \ref{lem:gamma2} for the cone $\overline{\ago_{0'}^{+}}$. 
The case of Proposition \ref{prop:decOf1Rel} is then in the same relation 
to the classical case as Lemma \ref{lem:gamma3} is to Lemma \ref{lem:gamma2}. 
The decomposition of $\overline{\ago_{0'}^{+}}$ is being given by the cones $\brzgop$.
The analogy is really close, Lemma \ref{lem:gamma3} takes Lemma \ref{lem:gamma2}
as a starting point and then refines it slightly, this what the proof of Proposition \ref{prop:decOf1Rel} 
essentially does too (although, we rather prove both at the same time). 

\begin{flushleft}
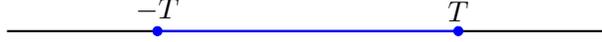
\begin{figure}
\begin{tikzpicture}[scale=2]
\draw[thick]  (-2, 0) -- (2, 0);
\draw[thick, blue]  (-1, 0) -- (1, 0);
\filldraw[blue] (1,0) circle (0.03cm);
\draw (1,0) node[above] {$T$};
\filldraw[blue] (-1,0) circle (0.03cm);
\draw (-1,0) node[above] {$-T$};
\end{tikzpicture}
  \caption{\tiny{The decomposition of a line that is used in \cite{ichYamGl, leMoi2}.}
  }
    \label{figure:9}
\end{figure}
\end{flushleft}

There's only one case in the literature, that we are aware of, that a non-trivial relative decomposition 
of $1$ was considered. It is in the work of Ichino-Yamana \cite{ichYamGl}, Lemma 2.3 and of the author \cite{leMoi2} Lemme 2.2 
(both are essentially the same)
in the case $G'= GL_n \sbs GL_{n+1} = G$ (diagonal inclusion).
Our construction is in fact different from these. The decomposition of loc. cit. 
takes the standard (Arthur's) decomposition on the big group $G$ and restricts it to $G'$. 
As mentioned above, we take the decomposition at the level of $G'$ as the starting point. Let us look at the cases $n=1$ and $n=2$ 
of $GL_n \sbs GL_{n+1}$ more closely. 
\begin{enumerate}
\item  If $n=1$, the Figure \ref{figure:6} represents the decomposition of Proposition \ref{prop:decOf1Rel}. 
Note that $F^{G}(x, T)$ is the characteristic function of $x \in GL_1(\A)$ such that $\log |x|_{\A} = T$ (up to normalization).
On the other hand, \cite{ichYamGl, leMoi2} induce the decomposition presented in figure \ref{figure:9}.
\item The case $n=2$ is discussed in \ref{ssec:gamma} and depicted in Figure \ref{figure:7}.
\end{enumerate} 

\bdem
Let $Q'= L'V'$ be the Levi decomposition of $Q'$. 
Using Iwasawa decomposition $G'(\A) = V'(\A)L'(\A)K'$
and writing $x = nmk$ accordingly, we see 
that the sum in question only depends on $m$. 
This shows that it is enough to consider the case $Q = G$, 
as the inclusion $M_{Q'} \hrar M_{Q}$ 
is of the same nature as $G' \hrar G$.

Let $x \in G'(\A)$ 
and let $\delta \in G'(F)$ be such that 
$\delta x \in \Sgl(G')$. 
In particular
\[
[\ago_{0'}^{+}](H_{0'}(\delta x) - \tzer) = 1.
\]

We apply part Lemma \ref{lem:gamma3}
to the cone $\overline{\ago_{0'}^{+}}$
and its decomposition 
\[
\overline{\ago_{0'}^{+}} = \bigcup \, \brzgop
\]
where the union is taken over $P \in \calP^{G}(P_{0}')$ such that $\zgo_{P}^{+}$ is open in $\ago_{0'}$,
to obtain
\[
[\ago_{0'}^{+}](H_{0'}(\delta x) - \tzer) = 
\sum_{P \in \calF^{G}(P_{0}')}\tlGam_{0'}^{P \cap G'}(H_{0'}(\delta x)^{P} - \tzer^{P}, T^{P} - \tzer^{P})
\tau_{P}(H_{0'}(\delta x)_{P} - T_{P}).
\]
The form of this decomposition is the real reason for the peculiar condition \eqref{eq:FPrel}. 
We obtain therefore that the sum in Proposition equals at least $1$.

Suppose now there exist $P, Q \in \calF^{G}(P_{0}')$ 
and $\delta_1 \in P'(F) \bsl G'(F)$
and  $\delta_2 \in Q'(F) \bsl G'(F)$, where $P' = P \cap G'$ and $Q'= Q \cap G'$, such that
$\delta_1 x \in \Sgl(P')$, $\delta_2 x \in \Sgl(Q')$ and both 
\[
\tlGam_{0'}^{P'}(H_{0'}(\delta_{1} x)^{P} - \tzer^{P}, T^{P} - \tzer^{P})\tau_{P}(H_{0'}(\delta_{1} x)_{P} - T_{P})
\]
and 
\[
\tlGam_{0'}^{Q'}(H_{0'}(\delta_{2} x)^{Q} - \tzer^{Q}, T^{Q} - \tzer^{Q})\tau_{Q}(H_{0'}(\delta_2 x)_{Q} - T_{Q})
\]
equal $1$.

Since $T$ is regular, it follows from Lemma \ref{lem:oldGamma} and 
the identity \eqref{eq:gammaOldNew} that both
\[
\tlGam_{0'}^{P'}(H_{0'}(\delta_{1} x)^{P} - \tzer^{P}, T^{P} - \tzer^{P})\tau_{P}(H_{0'}(\delta_{1} x)_{P} - T_{P}) 
\]
and 
\[
\tlGam_{0'}^{P'}(H_{0'}(\delta_{1} x)^{P'} - \tzer^{P'}, T^{P'} - \tzer^{P'})\tau_{P}(H_{P'}(\delta_{1} x) - T_{P'})
\]
equal $1$.
We obtain the same result for $Q$. Moreover, since 
$\zgo_{P}^{+} \sbs \ago_{P'}^{+}$ 
we obtain 
\begin{multline*}
\tlGam_{0'}^{P'}(H_{0'}(\delta_{1} x)^{P'} - \tzer^{P'}, T^{P'} - \tzer^{P'})[\ago_{P'}^{+}](H_{P'}(\delta_{1} x) - T_{P'})  = \\
\tlGam_{0'}^{Q'}(H_{0'}(\delta_{2} x)^{Q'} - \tzer^{Q'}, T^{Q'} - \tzer^{Q'})[\ago_{Q'}^{+}](H_{Q'}(\delta_2 x) - T_{Q'}) = 1.
\end{multline*}

Using Lemmas \ref{lem:oldGamma} and  
\ref{lem:labWal127} and \ref{lem:lang212} 
we get 
that $\delta_{2}\delta_{1}^{-1} \in P'(F) $ and $\delta_{1}\delta_{2}^{-1} \in Q'(F) $
which means that $\delta_{1}\delta_{2}^{-1} \in (P' \cap Q')(F) $. 
Let $R' = P' \cap Q'$. Let $H$ be the projection of $H_{0'}(\delta_1 x)$ 
onto $\ago_{R'}$. It equals the projection of $H_{0'}(\delta_2 x)$ onto $\ago_{R'}$ 
by what we have just seen. Using part 4 of Lemma \ref{lem:basicRootProps}
we obtain the equality
\[
\tlGam_{R'}^{P'}(H^{P'} - \tzer^{P'}_{R'}, T_{R'}^{P'} - \tzer^{P'}_{R'})
[\ago_{P'}^{+}](H_{P'} - T_{P'}) = \tlGam_{R'}^{Q'}(H^{Q'} - \tzer^{Q'}_{R'}, T_{R'}^{Q'} - \tzer^{Q'}_{R'})
[\ago_{Q'}^{+}](H_{Q'} - T_{Q'})
\]
The equality is between non-zero numbers and Lemmas \ref{lem:gamma2} and  \ref{lem:oldGamma} show that $P'$ must be $Q'$ for this to hold. 
Moreover, we get that $\delta_1$ and $\delta_2$ are equal mod $P'(F) = Q'(F)$ which proves the desired result. 
\edem

\bcor\label{cor:FPisSimple}
Let $x \in G'(\A)$ and let $\delta \in Q'(F)$ 
be such that $\delta x \in \Sgl^{G'}(Q')$. 
Then 
\[
F^{Q}(x, T) = \tlGam_{0'}^{Q'}(H_{0'}(\delta x)^{Q} - \tzer^{Q}, T^{Q} - \tzer^{Q}).
\]
In particular, the right hand side doesn't depend on the choice of $\delta$.
\ecor
\bdem
If $F^{Q}(x, T) = 0$ the result is clear. 
If there exist  $\delta, \delta' \in Q'(F)$ 
such that $\delta x, \delta' x \in \Sgl(Q') $ 
and 
\[
 1 =  \tlGam_{0'}^{Q'}(H_{0'}(\delta x)^{Q} - \tzer^{Q}, T^{Q} - \tzer^{Q})
 \neq  \tlGam_{0'}^{Q'}(H_{0'}(\delta' x)^{Q} - \tzer^{Q}, T^{Q} - \tzer^{Q}).
\]
we obtain a contradiction with Proposition \ref{prop:decOf1Rel}
\edem

We also note the following immediate consequence. 

\bcor\label{cor:fqComp} Let $x \in G'(\A)$ be 
such that $F^{Q}(x,T) = 1$. Then, there exists a $\gamma \in Q'(F)$ 
such that $\gamma x$ can be written as $nac$ 
where $n \in N_{Q'}(\A)$, $a \in Z_{Q}^{\infty}$ 
and $c$ belongs to a compact in $G'(\A)$ that depends only on $T$.
\ecor

\subsection{Relative truncation operator}\label{ssec:relTrunc}

Let $Q \in \calF^{G}(P_{0}')$, 
and let $Q' = Q \cap G$.
Let $\phi$ be a locally integrable function on $Q(F) \bsl G(\A)$. 
We define the relative (mixed) truncation operator as 
\[
\La^{T,Q}\phi(x) = 
\ssum{
P \in \calF^{G}(P_{0}') \\
P \sbs Q 
}
\varepsilon_{P}^{Q} \sum_{\delta \in (G' \cap P)(F) \bsl Q'(F)}
\htau_{P}^{Q}(H_{0'}(\delta x)^{Q} - T^{Q})\phi_{P}(\delta x), \quad x \in Q'(F) \bsl G'(\A).
\]
When $Q = G$ we write simply $\La^{T} = \La^{T,G}$. Note that the sums in the
 the definition of the operator $\La^{T,Q}$ are all finite 
thanks to Lemma \ref{lem:htauFinite}.

We first establish some classical combinatorial properties of our operators. 

\blem\label{lem:invForm} For all locally integrable function $\phi$ on $[G]$ 
and all $x \in G'(\A)$ 
we have the identity
\[
\phi(x) = \sum_{Q \in \calF^{G}(P_{0}')}\sum_{\delta \in (Q \cap G')(F) \bsl G'(F)}
\tau_{Q}(H_{0'}(\delta x)_{Q}  - T_{Q})\La^{T,Q}\phi(\delta x).
\]
\elem

\bdem
We have
\begin{multline*}
\sum_{Q \in \calF^{G}(P_{0}')}\sum_{\delta \in (Q \cap G')(F) \bsl G'(F)}
\tau_{Q}(H_{0'}(\delta x)_{Q} - T_{Q})\La^{T,Q}\phi(\delta x) = \\
\sum_{P \in \calF^{G}(P_{0}')}
\sum_{\delta \in (P \cap G')(F) \bsl G'(F)}\phi_{P}(\delta x)
\ssum{
Q \in \calF^{G}(P_{0}') \\
Q \sps P 
}
\varepsilon_{P}^{Q} 
\tau_{Q}(H_{0'}(\delta x)_{Q}  - T_{Q})
\htau_{P}^{Q}(H_{0'}(\delta x)^{Q} - T^{Q})
\end{multline*}
By part 1. of Proposition \ref{prop:langLemma} the last sum above is zero unless 
$P=G$. 
\edem

\blem\label{lem:laTT} With assumptions as above, for all $T,T'  \in \ago_{0'}$
we have the identity
\[
\La^{T+T'}\phi(x) = \sum_{Q \in \calF^{G}(P_{0}')}\sum_{\delta \in (Q \cap G')(F) \bsl G'(F)}
\Gamma_{Q}(H_{0'}(\delta x)^{G}_{Q} - (T')^{G}_{Q}, T^{G}_{Q})\La^{T',Q}\phi(\delta x).
\]
\elem
\bdem
This follows immediately from part 3. of Lemma \ref{lem:conesProps}.
\edem

The next theorem is the fruit of all constructions and results up to this point.

\btheo\label{thm:rapDec} 
Let $Q \in \calF^{G}(P_{0}')$ and set $Q' = Q \cap G'$.
Let $\phi$ be a smooth function on $Q(F) \bsl G(\A)$ 
of uniform moderate growth. 
Then, for all $T \in \treg + \ago_{0'}^{+}$ and $N > 0$ 
there exists a constant $C$ such that 
for all $x \in \Sgl(M_{Q'}) \cap M_{Q'}(\A)^{1, Q}$ and all $k \in K'$ we have 
\[
|\La^{T, Q} \phi(xk)| \le C |x|^{-N}.
\]
\etheo
\bdem
It is enough to take $Q = G$ and $k = 1$.

Let $x \in G'(\A)^{G,1}$. Then 
$H_{0'}(x)^{G} = H_{0'}(x)$. We can also assume $T = T^{G}$ as $\La^{T} = \La^{T^{G}}$.
Using Proposition \ref{prop:decOf1Rel} and the relation (4) of Lemma \ref{lem:conesProps} 
we obtain
\begin{multline*}
\La^{T}\phi(x) = 
\sum_{R \in \calF^{G}(P_{0}')}
\varepsilon_{R}^{G}\sum_{\delta \in R'(F) \bsl G'(F)}
\htau_{R}(H_{0'}(\delta x) - T)\phi_{R}(\delta x) = 
\sum_{R \in \calF^{G}(P_{0}')}
\varepsilon_{R}^{G}\sum_{\delta \in R'(F) \bsl G'(F)} \\
\dsl \ssum{
P \in \calF^{G}(P_{0}')\\
P \sbs R
} \sum_{\xi \in P'(F) \bsl R'(F)}F^{P}(\xi \delta x, T)\tau_{P}^{R}(H_{0'}(\xi \delta x)_{P} - T_{P})
\rb 
\htau_{R}(H_{0'}(\delta x) - T)\phi_{R}(\delta x) = \\
\ssum{
P, R \in \calF^{G}(P_{0}')\\
P \sbs R
} \varepsilon_{R}^{G}
\sum_{\delta \in P'(F) \bsl G'(F)}
F^{P}(\delta x, T)\tau_{P}^{R}(H_{P}(\delta x) - T_{P})
\htau_{R}(H_{0'}(\delta x)_{P} - T_{P})\phi_{R}(\delta x) = \\
\ssum{
P, Q, R \in \calF^{G}(P_{0}')\\
P \sbs R \sbs Q
} \varepsilon_{R}^{G}
\sum_{\delta \in P'(F) \bsl G'(F)}
F^{P}(\delta x, T)\sigma_{P}^{Q}(H_{0'}(\delta x)_{P} - T_{P})\phi_{R}(\delta x) = \\
\ssum{
P, Q\in \calF^{G}(P_{0}')\\
P \sbs Q
} \varepsilon_{Q}^{G}
\sum_{\delta \in P'(F) \bsl G'(F)}
F^{P}(\delta x, T)\sigma_{P}^{Q}(H_{0'}(\delta x)_{P} - T_{P})
\dsl
\ssum{
R \in \calF^{G}(P_{0}')\\
P \sbs R \sbs Q
} \varepsilon_{R}^{Q}\phi_{R}(\delta x)
\rb
\end{multline*}
where $P' = P \cap G'$, $R'= R \cap G'$ 
and we used the fact that for $P \sbs R$
\[
\htau_{R}(H_{P}) = \htau_{R}(H), \quad H \in \ago_{0'}.
\]

It is enough to prove
\[
\sum_{\delta \in P'(F) \bsl G'(F)}
F^{P}(\delta x, T)\sigma_{P}^{Q}(H_{0'}(\delta x)_{P} - T_{P})
\dsla
\ssum{
P \sbs R \sbs Q\\
R \in \calF^{G}(P_{0}')
} \varepsilon_{R}^{Q}\phi_{R}(\delta x)
\rba 
\le C |x|^{-N}
\]
for all  $P, Q \in \calF^{G}(P_{0}')$ such that $P \sbs Q$ and some $C > 0$.

Suppose first $P = Q = G$. 
Then, the inequality 
\[
|F^{G}(x,T)\phi(x)| \le C |x|^{-N}
\]
follows from Corrollary \ref{cor:FPisSimple}. 
If $P = Q \neq G$, we have $\sigma_{P}^{Q} \equiv 0$ 
thanks to Proposition \ref{prop:langLemma} part 2. 

Assume from now on then that $P \sbn  Q$ 
and consider $P$ and $Q$ fixed. 
Let $R \in  \calF^{G}(A_{0})$ 
be such that $R \sbs Q$. 
Define $R^{\sharp}$ to be the smallest  element of $\calF^{G}(P_{0}')$ contained in $Q$ 
containing $R$. It is a well defined notion 
since for any two $P_1, P_2 \in \calF^{G}(P_{0}')$ 
we have 
\[
P_1 \cap P_2 \in \calF^{G}(P_{0}) \ \Longrightarrow \ P_1 \cap P_2 \in \calF^{G}(P_{0}'). 
\]

Applying the notation of Pararagraph \ref{ssec:altSum} and the identity \eqref{eq:constAlt} therein, 
we obtain
\begin{multline*}
\ssum{
R \in \calF^{G}(P_{0}') \\
P \sbs R \sbs Q
} \varepsilon_{R}^{Q}\phi_{R}= 
\ssum{
R \in \calF^{G}(P_{0}') \\
R_{1} \in \calF^{G}(A_{0}) \\
P \sbs R_1 \sbs R \sbs Q
} \varepsilon_{R}^{Q}\phi_{P, R_1} = \\
\ssum{
R_{1} \in \calF^{G}(A_{0}) \\
P \sbs R_1
}\phi_{P, R_1} 
\ssum{
R \in \calF^{G}(P_{0}')\\
R_1 \sbs R \sbs Q 
} \varepsilon_{R}^{Q} = 
\ssum{
R_{1} \in \calF^{G}(A_{0}) \\
P \sbs R_1}\phi_{P, R_1} 
\ssum{
R \in \calF^{G}(P_{0}')\\
R_1^{\sharp} \sbs R \sbs Q 
} \varepsilon_{R}^{Q}.
\end{multline*}
It follows from \eqref{eq:euler} that it is enough to show
\[
\sum_{\delta \in P'(F) \bsl G'(F)}
F^{P}(\delta x, T)\sigma_{P}^{Q}(H_{0'}(\delta x)_{P} - T_{P})
|\phi_{P, R}(\delta x)|
\le C |x|^{-N}
\]
for all $R \in \calF^{G}(A_{0})$ 
such that $P \sbs  R\sbs Q$ and $R^{\sharp} = Q$.

Let $\delta \in P'(F) \bsl G'(F)$ be such that $F^{P}(\delta x, T)\sigma_{P}^{Q}(H_{0'}(\delta x)_{P} - T_{P}) = 1$. 
Using Corollary \ref{cor:fqComp} we can write 
\[
\delta x = n' a c
\]
with $n' \in N'(\A)$, $a \in Z_{P}^{\infty}$ 
and $c$ in a fixed compact of $G'(\A)$. Multiplication on the right by elements belonging 
to a fixed compact only influence the multiplicative constant in the estimates 
we can thus assume $c = 1$.

Using Lemma \ref{lem:431} we have that for all $r \in \N$ 
there exists a constant $c_1$ and a finite number of differential operators $X_k$ 
such that  
\begin{equation}\label{eq:lem431}
|\phi_{P,R}(n'a)| \le \sup_{n \in [N_{P}]}|\phi_{P,R}(na)| \le c_{1}\sum_{k}\sup_{n \in [N_{P}]}
|\rho(\Ad(a^{-1})X_{k}) \phi(na)| 
\end{equation}
where for all $k$ there exists a $\la_{k} = \sum_{\al \in \Delta_{P}^{R}}a_{k,\al}\al$
with $a_{k,\al} \in \N$ such that $a_{k,\al} \ge r$ 
and
\[
\Ad(a^{-1})X_{k} = 
e^{-\langle \la_{k},  H_{P}(a) \rangle }
X_{k}.
\]
Note that since $a \in Z_{P}^{\infty}$ we actually have $H_{P}(a) = H_{0'}(a)_{P} \in \zgo_{P}$. 
Let $H = H_{0'}(a)_{P}$ and let $H_1 + H_2$ be its decomposition along $\zgo_{P} = \zgo_{P}^{Q} \oplus \zgo_{Q}$.

We have $\Delta_{P}^{R} \sbs \Delta_{P}^{Q}$. 
Therefore $\langle \la,  H_{P}(a) \rangle  = \langle \la,  H_1 \rangle $. 
The elements of $\Delta_{P}^{Q}$ 
are (strictly) positive on $\rint A(\bragoq, \bragop)$ (following the description \eqref{eq:angRoot})
which contains $\rint A(\brzgoq, \brzgop)$.
Moreover, since $R^{\sharp} = Q$, the restriction of elements of $\Delta_{P}^{R}$ to $\zgo_{P}^{Q}$ 
spans the linear dual of the latter. Therefore, using the first assertion of 
part 5 of Lemma \ref{lem:conesProps}, 
there exists a constant $c_2$ such that 
\[
e^{-\langle \la,  H_1 \rangle} \le e^{-rc_2\|H_1\|}.
\]
The second assertion of Lemma \ref{lem:conesProps} part 5, shows that
\[
e^{-rc_2\|H_1\|} \le e^{-rc_3\|H\|}
\]
 for some constant $c_3$. 
Splitting this expression in two and using the uniform moderate growth hypothesis on $\phi$
we obtain
\[
F^{P}(\delta x, T)\sigma_{P}^{Q}(H_{0'}(\delta x)_{P} - T_{P})
|\phi_{P, R}(\delta x)|
\le c_4 |\delta x|^{-N} \le c_5 |x|^{-N} 
\]
for any fixed $N$ and some constants $c_4$, $c_5$ that depend on $N$ and $T$.

We still have to take the sum over $\delta \in P'(F) \bsl G'(F)$. 
We have
\[
|\sigma_{P}^{Q}(H_{0'}(\delta x)_{P} - T_{P})| \le 
\ssum{
S \in \calF^{G}(P_{0}') \\
S \sps Q
}
\tau_{P}^{S}(H_{0'}(\delta x)_{P} - T_{P})\htau_{S}(H_{0'}(\delta x)_{S} - T_{S}) .
\]
Fix then $S \in \calF^{G}(P_{0}')$ containing $Q$. 
We want to show that the number 
of $\delta \in P'(F) \bsl G'(F)$ 
such that 
\[
F^{P}(\delta x, T)
 \tau_{P}^{S}(H_{0'}(\delta x)_{P} - T_{P}) \htau_{S}(H_{0'}(\delta x)_{S} - T_{S}) \neq 0
\]
is bounded by a power of $|x|$, independent of $T$. 

Let $S'= S \cap G'$. 
We have
\begin{multline*}
\sum_{\delta \in P'(F) \bsl G'(F)}
F^{P}(\delta x, T) \tau_{P}^{S}(H_{0'}(\delta x)_{P} - T_{P}) \htau_{S}(H_{0'}(\delta x)_{S} - T_{S}) = \\
\sum_{\gamma \in S'(F) \bsl G'(F) }
\sum_{\xi \in P'(F) \bsl S'(F)}
F^{P}(\xi \gamma x, T)\tau_{P}^{S}(H_{0'}(\xi \gamma x)_{P} - T_{P})\htau_{S}(H_{0'}(\xi \gamma x)_{S} - T_{S}).
\end{multline*}
By Proposition \ref{prop:decOf1Rel}, for each $\gamma \in S'(F) \bsl G'(F) $, there 
is at most one $\xi \in P'(F) \bsl S'(F)$ 
such that 
\[
F^{P}(\xi \gamma x, T)\tau_{P}^{S}(H_{0'}(\xi \gamma x)_{P} - T_{P}) \neq 0.
\]
It is enough therefore to consider 
\[
\sum_{\gamma \in S'(F) \bsl G'(F) }
\htau_{S}(H_{0'}(\gamma x)_{S} - T_{S}).
\]
This sum is bounded by a power of $|x|$ independent of $T$ in virtue of Lemma \ref{lem:htauFinite}, 
which ends the proof. 
\edem

\section{Invariant period}\label{sec:invPeriod}
At this point we have proven all the results pertaining to convergence 
as well as the various combinatorial result about truncation. 
We are ready to define regularization of the period integrals of automorphic forms. 
Everything has been set up so that we can use the formalism of \cite{jlr} 
almost verbatim. Therefore, most of the proofs here will be brief.

\subsection{The $\rho_{P}$}
Recall that for all $P \in \calF^{G}(A_0)$ we have an element $\rho_{P} \in \Hom_{\R}(\ago_{0}, \R) \cong \ago_{0}$ 
defined as the half sum of roots for the action of $A_{P}$ on $\Lie (N_{P})$.

\subsection{Automorphic forms}

Let $P \in \calF^{G}(A_0)$. Let $\calA_{P}(G)$ be the set of automorphic forms on $N(\A)M(F) \bsl G(\A)$. 
These are smooth, moderate growth functions $\phi : N(\A)M(F) \bsl G(\A) \to \C$ that are $K$-finite 
and finite for the action of the center of $\Ugo(G)$. Automorphic forms are automatically of uniform moderate growth.
If $P=G$ we write $\calA(G)$ for $\calA_{G}(G)$.

Every automorphic form $\phi \in \calA_{P}(G)$ admits a fine sum decomposition
satisfying
\[
\phi(ay) = \sum_{i} q_{i}(H_{P}(a))e^{ \langle \la_{i} + \rho_{P},H_{P}(a) \rangle}\phi_{i}(y)
\]
for $a \in A_{P}^{\infty}$ and $y \in N(\A)M(\A)^{1}K$, 
with $q_{i} \in \C[\ago_{P}]$, $\la_{i} \in \ago_{0,\C}$ and $\phi_{i} \in \calA_{P}(G)$
such that $\phi_{i}(ax) = \phi_{i}(x)$ for $a \in A_{P}^{\infty}$ and all $x \in G(\A)$.
The set composed of distinct $\la_{i}$ is uniquely determined by $\phi$ and is called the set of exponents of $\phi$. 
For $Q \sbs P$ the set of exponents of $\phi$ along $Q$ is by definition 
the set of exponents of $\phi_{Q}$. It will be denoted by $\calE_{Q}(\phi)$.
 
Let $P \in \calF^{G}(P_{0}')$ and $P'  = P \cap G$.
Let $\phi \in \calA(G)$. 
In the relative setting only restrictions of exponents to $\ago_{0',\C}$ matter, so let  us denote 
\[
\calE_{P}(\phi)' \sbs \ago_{0',\C}
\]
the set of restrictions (projections) of elements of $\calE_{P}(\phi)$ to $\ago_{0',\C}$. 
It follows that for every $\phi \in \calA(G)$, its constant term 
with respect to $P$ admits a decomposition satisfying 
\begin{equation}\label{eq:autDecAdpt}
\phi_{P}(ay) = \sum_{i} q_{i}(H_{0'}(a)_{P})e^{\langle \la_{i} + \rho_{P}, H_{0'}(a)_{P} \rangle}\phi_{i}(y)
\end{equation}
for $a \in Z_{P}^{\infty}$ and $y \in N'(\A)M'(\A)^{1, P}K'$, 
with $q_{i} \in \C[\zgo_{P}]$, $\la_{i} \in \calE_{P}(\phi)'$ and $\phi_{i} \in \calA_{P}(G)$ 
that are left $Z_{P}^{\infty}$-invariant.

\subsection{Haar measures}

We fix a Haar measure on $G'(\A)$ and on $K'$ giving it the volume $1$. 
For all unipotent subgroups $N'$ of $G'$ we fix the Haar measure on 
$N'(\A)$ giving $[N']$ measure $1$. 
We have a natural measure on $\ago_{0'}$ coming from the Euclidean structure, 
it induces measures on all subspaces of $\ago_{0'}$. 
In particular, for all $P \in \calF^{G}(P'_0)$, the isomorphism $H_P : Z_{P}^{\infty} \to \zgo_{P}$ 
induces a measure on $Z_{P}^{\infty}$. If we note $M' = M \cap G'$, it follows that there is a unique 
Haar measure on $M'(\A)^{1,P}$ such that for all $f \in L^{1}(G'(\A))$
\[
\int_{G'(\A)}f(g)\,dg = 
\int_{K'} \int_{N'(\A)}\int_{M'(\A)^{1,P}}\int_{Z_{P}^{\infty}}
e^{\langle -2\rho_{P'}, H_{P'}(a) + H_{P'}(m)\rangle}f(nmak)\,
dadmdndk.
\]

\subsection{Polynomial exponentials}

 We will use the results and notation of Section \ref{ssec:FT}.
Let $\calF^{G, max}(P_{0}')$ be the set of $P \in \calF^{G}(P_{0}')$ 
such that $\dim \zgo_{P}^{G} = 1$.

Let $P \in \calF^{G}(A_0)$ and $q \in \C[\zgo_{P}]$. 
We have the meromorphic function 
\begin{equation*}
\scrF(\brzgop, q, \la) \quad \la \in \ago_{0',\C}.
\end{equation*}
The function $\la \to \scrF(\brzgop, q, \la)$ is defined by a convergent integral for 
\[
\Rel(\la) \in -\rint (\zgo_{P}^{G} \cap \brzgop)^{\vee}.
\]
Moreover, the function $\scrF(\brzgo_{P}^{+}, q,\la)$ is holomorphic 
on the open set
\[
\ago_{0', \C}^{P-reg} := \{\la \in \ago_{0',\C} \ | \ \langle \la, \zgo_{Q}^{G} \rangle \neq 0, \  
\forall \, Q \in \calF^{G, max}(P_{0}') \cap \calF^{G}(P) \}.
\]

For all $T \in \ago_{0'}$, we also have the holomorphic function 
\begin{equation*}
\scrF(\Gamma_{P}, T, q, \la) := \int_{\zgo_{P}^{G}}\Gamma_{P}(H, T)e^{\langle \la, H \rangle}q(H)\,dH \quad \la \in \ago_{0',\C}.
\end{equation*}
As a function of the variable $T$, $\scrF(\Gamma_{P}, T, q, \la)$ is a polynomial-exponential 
whose purely polynomial term is constant for $\la \in \ago_{0', \C}^{P-reg}$ 
and given by
\[
\scrF(\brzgop, q, \la).
\]

\subsection{The period}

Let $\xi : G'(\A) \to \C^{\times}$ be an automorphic ($G'(F)$-invariant) 
character. 
By abuse of notation, we will also denote by $\xi$, the unique
element of $\ago_{G',\C}$
such that 
\[
\xi(a) = e^{\langle \xi, H_{G'}(a) \rangle} \quad a \in A_{G'}^{\infty}.
\]
 
For $P \in \calF^{G}(P_{0}')$ let 
\[
\upla_{P}  \in \ago_{0'}
\] 
be the projection of $\rho_{P}- 2\rho_{P \cap G'}$ onto $\ago_{0'}$.

Let $\calA(G)^{\xi-reg}$ be the space of $\phi \in \calA(G)$ such that 
for all $P \in \calF^{G}(P_{0}')$ and all $\la \in \calE_{P}(\phi)'$ 
we have
\[
\la + \xi + \upla_{P} \in \ago_{0', \C}^{P-reg}.
\]
 In fact, $\calA(G)^{\xi-reg}$ is simply the space of 
$\phi \in \calA(G)$ such that 
for all $P \in \calF^{G, max}(P_{0}')$ and all $\la \in \calE_{P}(\phi)'$ 
we have
\[
\la + \xi + \upla_{P} \in \ago_{0', \C}^{P-reg} \Leftrightarrow 
\langle \la + \xi + \upla_{P} , \zgo_{P}^{G} \rangle \neq 0.
\]
Note that $\calA(G)^{\xi-reg}$ is $G(\A)$-stable.

For $T \in \ago_{0'}$ we define the truncated period to be the following functional on $\calA(G)$
\[
\calP^{T}(\xi, \phi) := \int_{[G']^{1,G}}\La^{T} \phi(h)\xi(h)\,dh.
\]
In general, for $P \in \calF^{G}(P_{0}')$ and $\phi \in \calA_{P}(G)$ set
\[
\calP^{P, T}(\xi, \phi) := \int_{K'}\int_{[M']^{1,P}}e^{\langle -2\rho_{P'}, H_{P'}(m) \rangle}\La^{T,P}\phi(mk)\xi(mk)\,dmdk
\]
where $P'= G' \cap P$.
Theorem \ref{thm:rapDec} asserts that $\calP^{P, T}(\xi, \cdot)$ is well defined for $T \in \treg+ \ago_{0'}^{+}$.

\btheo\label{theo:theProp} The following assertions hold
\begin{enumerate}
\item There exists a unique polynomial-exponential function on $\ago_{0'}$ 
that coincides with $T \mapsto \calP^{T}(\xi, \phi)$ for $T \in \treg + \ago_{0'}^{+}$. 
\item For $\phi \in \calA(G)^{\xi-reg}$ the purely polynomial part of $\calP^{T}(\xi, \phi)$ 
is constant. 
\item  
Define the functional $\calP(\xi, \cdot)$ on $\calA(G)^{\xi-reg}$ 
as the constant term of $T \to \calP^{T}(\xi, \cdot)$.
Then, for all $\phi \in \calA(G)^{\xi-reg}$ and all $T'\in \treg + \ago_{0'}^{+}$ we have
\begin{equation}\label{eq:invPer}
\calP(\xi, \phi) = \sum_{P}\sum_{i}
e^{\langle \la_{i} + \upla_{P} + \xi, T'_{P} \rangle}
\scrF(\brzgo_{P}^{+}, (q_{i})_{T'}, \la_{i} + \upla_{P} + \xi)
\calP^{P, T'}(\xi, \phi_{i})
\end{equation}
where we write $\phi_{P} \in \calA_{P}(G)$ as in \eqref{eq:autDecAdpt}
and $(q_{i})_{T'} (H) = q_{i}(H+(T')^{G}_{P})$. In particular, 
the sum in \eqref{eq:invPer} is independent of $T'$.
\item The functional $\calP(\xi, \cdot)$ is right $G'(\A)$-$\xi$-equivariant, that is
\[
\calP(\xi, \phi_{x}) = \xi(x)\calP(\xi, \phi), \quad \phi \in \calA(G)^{\xi-reg}, \ x \in G'(\A)
\]
where $\phi_{x}(y) = \phi(yx^{-1})$. Moreover, it is independent of 
the various choices made in the construction of the operator  $\La^{T}$.
\end{enumerate}
\etheo
\bdem
Let's prove the points 1,2 and 3. Fix $T' \in \ago_{0'}^{G}$. We have then, using Lemma \ref{lem:laTT}  
\[
\calP^{G, T + T'}(\xi, \phi) = \sum_{P \in \calF^{G}(P_{0}')}
\sum_{i}e^{\langle \la_{i} + \upla_{P} + \xi, T'_{P} \rangle}
\scrF(\Gamma_{P}, T^{G}_{P}, (q_{i})_{T'}, \la_{i} + \upla_{P} + \xi)\calP^{P, T'}(\xi, \phi_{i})
\]
where we write $\phi_{P}$
as in \eqref{eq:autDecAdpt} and $(q_{i})_{T'} (H) = q_{i}(H+T'_{P})$.
The statement follows now from the properties of the Fourier transform of $\Gamma_{P}$. 
In particular, we see that the purely polynomial part 
of $T \to \calP^{G, T + T'}(\xi, \phi)$
is constant 
when $\la_{i} + \xi + \upla_{P} \in \ago_{0', \C}^{P-reg}$ for all $P \in \calF^{G}(P_{0}')$
and given by the formula in point 3. 

Let's prove part 4. 
For all $y \in G'(\A)$ let $K(y) \in K'$ be any element 
such that $yK(y)^{-1} \in P_{0}'(\A)$. 
Fix $x \in G'(\A)$. 
Just as in \cite{jlr}, equation (26) we have for all $y \in G'(\A)$
\[
\La^{T}(\phi_{x^{-1}})(yx) = \sum_{P \in \calF^{G}(P_{0}')}
\sum_{\delta \in (P \cap G')(F) \bsl G'(F)}
\Gamma_{P}(H_{0'}(\delta y)_{P} - T_{P}, -H_{0'}(K(\delta y)x)_{P})
\La^{T,P}\phi(\delta y).
\]
It follows immediately after integrating both sides 
that constant terms (i.e. purely polynomial parts in $T$) of $\calP^{G, T}(\xi, \phi_{x^{-1}})$
and $\xi(x)\calP^{G, T}(\xi, \phi)$ coincide. 

The proof of independence of various choices relies on the same principle and will be omitted. 
\edem

\bcor Suppose $\zgo_{G} \neq \ago_{G'}$. Let $\phi \in \calA(G)$.
Then $\phi \in \calA(G)^{\xi-reg}$ for almost all $\xi \in \ago_{G',\C}$. 
Moreover, $\xi \in \ago_{G',\C} \to \calP(\xi, \phi)$ is a meromorphic function with hyperplane singularities on $\ago_{G',\C}$.
\ecor
\bdem 
Let $\ago_{G'}^{G}$ be the orthogonal complement of $\zgo_{G}$ in $\ago_{G', \C}$.
It suffices to realize that generic elements of $\ago_{G'}^{G}$ are not zero on $\zgo_{P}^{G}$ 
for all $P \in \calF^{G, max}(P_{0}')$.
\edem

\subsection{The case of non-connected $G'$}\label{ssec:nonConn}

Our construction also applies to non-connected groups. 
Let $\tlG' \sbs G$ be a reductive algebraic subgroup of $G$. Let $G'$ be the connected component of $\tlG'$. 
Let $\calP_{G'}$ be the regularized period as defined in Proposition \ref{theo:theProp} with respect to $G' \sbs G$ 
and a trivial character $\xi$ (for simplicity). Let $\calA(G)^{reg}$ be the associated subspace of definition of $\calP_{G'}$

For $\phi \in \calA(G)$ and $h \in G(\A)$ let $\phi_{h} \in \calA(G)$ be the right translation of $\phi$ by $h^{-1}$. 
Note that, by Tychonoff's theorem, the quotient $G'(\A) \bsl \tlG'(\A)$ is compact. 
We see then that the following functional
\[
\calP_{\tlG'}(\phi) = (\# (G(F) \bsl \tlG (F)))^{-1} \int_{G(\A) \bsl \tlG (\A)} \calP(\phi_{h})\,dh , \quad \phi \in \calA(G)^{reg}
\]
is well defined and $\tlG'(\A)$-equivariant. 

\subsection{Periods of Eisenstein series induced from maximal parabolic subgroups}\label{ssec:cuspEis}

Since we will use the results of this paper in \cite{polWanZyd2} in the particular case of 
Eisenstein series induced from maximal parabolic subgroup of $G$ we include this case here. 
We hope this Paragraph will also demonstrate that the constructions 
of this article are very explicit in practice. 

Assume for simplicity that $A_{G}^{\infty}$ is trivial. 
Let $P \in \calF^{G}(A_{0})$ be a maximal, self-associate,
parabolic subgroup
of $G$. 
Let $\phi \in \calA_{P}(G)$ be an automorphic form 
such that $m \in M(\A) \to \phi(mg)$ 
is a cusp form and
$\phi(ag) = e^{\rho_{P}( H_{P}(a)) } \phi(g)$
for all $g \in G(\A)$ and $a \in A_{P}^{\infty}$.
Let $E(g, \phi, s) $ be the associated Eisenstein series 
where $s \in \C$ and we fix the isomorphism $\C \cong \ago_{P,\C}$ 
that identifies $1$ with the fundamental weight $\varpi_{P} \in \hDelta_{P}$. 
Let $c \in \Q_{+}^{\times}$ be such that
\[
\rho_{P} = c\varpi_{P}.
\]

Let $Q \in \calF^{G}(A_0)$ be a parabolic 
subgroup conjugate to $P$. 
We let $\Omega(P,Q)$ be the two element set of isometries
between $\ago_{P}^{G}$ and $\ago_{Q}^{G}$. 
To $w \in \Omega(P,Q)$ 
we assign its sign $\sgn(w) \in \{-1, 1\}$ 
such that $w\varpi_{P} = \sgn(w) \varpi_{Q}$.
We have then that 
\[
E(\phi, s\varpi, x)_{Q} = 
\sum_{w \in W(P,Q)}
M(w,s)\phi(x)e^{\langle  \sgn(w)s \varpi_{Q}, H_{Q}(x) \rangle} 
\]
where $M(w, s) : \calA_{P}(G) \to \calA_{Q}(G)$ is an intertwining operator, independent of $s$ if $\sgn(w) = 1$.

Let $Q \in \calF^{G}(A_0)$ be as above and 
suppose moreover $Q \in \calF^{G}(P_{0}')$. Set $Q'= Q \cap G'$.
Since $\ago_{G} = \{0\}$, 
the space $\zgo_{Q} = \zgo_{Q}^{G}$ is one dimensional and because it intersects $\ago_{Q}^{+}$ 
we must have $\zgo_{Q} = \ago_{Q}$. 
Let $c_{Q}^{G'}$ be such that $\rho_{Q'}$ projected onto $\zgo_{Q}$
equals $c_{Q}^{G'}\rho_{Q}$. 
For example, if $N_{Q}$ is abelian we have $c_{Q}^{G'} = \dfrac{\dim_{F}N_{Q'}}{\dim_{F}N_{Q}}$. 
For these data, the equation \eqref{eq:invPer} of Theorem \ref{theo:theProp} gives

\bcor\label{cor:cuspEis} The regularized period $\calP(E(\phi, s)) $ equals
\begin{multline*}
\int_{[G']^{1,G}}\La^{T}E(\phi, s\varpi, h)\,dh - \\
\sum_{Q}\sum_{w \in W(P,Q)}\dfrac{e^{\langle (\sgn(w)s + c(1-2c_{Q}^{G'}) )\varpi_{Q}, T \rangle}}{\sgn(w)s + c(1-2c_{Q}^{G'})}
\int_{K'}\int_{[M_{Q'}]^{1,Q}}e^{\langle -2\rho_{Q'}, H_{Q'}(m) \rangle}M(w,s)\phi(mk)\,dmdk
\end{multline*}
where the sum runs over $Q \in \calF^{G}(P_{0}')$ conjugate to $P$.
\ecor

\subsection{Critetirion for integrability}\label{ssec:intCrit}
We will need a bit of notation in this paragraph.
We fix a standard minimal parabolic subgroup $P_0 = M_{0}N_{0}$ of $G$, which we didn't need to fix thus far. 
We will use notation of of Paragraph \ref{ssec:altSum}. 
For any $P \in  \calF^{G}(P_{0})$ we set 
\[
\al_{P} = \sum_{\al \in \Delta_{0} \smin \Delta_{0}^{P}} \al \in \ago_{0}.
\]
We say  that $\phi : G(\A) \to \C$ has a central character 
if $\phi(ax) = e^{\langle \zeta , H_{G}(a) \rangle }\phi(x)$ for some $\zeta \in \ago_{G,\C}$ 
and all $a \in A_{G}(\A)$ and $x \in G(\A)$.

We start with invoking Lemma I.2.10 of \cite{moeWald2}. 

\blem\label{lem:I210}  
Let $P \in \calF^{G}(P_{0})$ be maximal. 
Fix $K_{f}$,  an open compact subgroup of $G(\A^{\infty})$, $X \in \Ugo(G)$ and $t > 0$. 
Then, there exists a finite collection $\{X_{i}\}_{i \in I} \sbs \Ugo(G)$
such that for all smooth functions $\phi$ on $N_{0}(F)\bsl G(\A) / K_{f}$ with a central character
the inequality
\begin{equation}\label{eq:I210}
|\rho(X_{i})\phi(x)| \le c_{i}e^{\langle \la, H_{0}(x) \rangle}, \quad \forall \, x \in \Sgl(G) \cap G(\A), \ \forall \, i \in I
\end{equation}
for some $\la \in ago_{0}$ and constants $c_{i}$ implies 
\[
|\rho(X)\phi_{P,G}(x)| \le c e^{\langle \la - t\al_{P},  H_{0}(x) \rangle}, \quad \forall \, x \in \Sgl(G) \cap G(\A)
\]
for some constant $c$. 
\elem

\bcor\label{cor:I210}  
Let $K_{f} \sbs G(\A^{\infty})$ be an open compact subgroup. 
Let $\phi$ be a function on $N_{0}(F)\bsl G(\A) / K_{f}$ of uniform moderate growth with a central character.
Then, there exists a $\mu \in \ago_{0}$
such that for all for all $X \in \Ugo(G)$, $t > 0$ 
and $P \in \calF^{G}(P_{0})$ 
there exists a $c >0$ 
such that 
\[
|\rho(X)\phi_{P,G}(x)| \le c e^{\langle \mu - t\al_{P},  H_{0}(x) \rangle}\quad \forall \, x \in \Sgl(G) \cap G(\A).
\]
\ecor
\bdem
The result is trivial if $P$ equals $G$ and  if $P$ is maximal then it follows from \ref{lem:I210} 
by uniform growth hypothesis.
We will prove the result by induction on the number of standard maximal parabolic subgroups containing $P$.
Let $P_1, \ldots, P_{k} \in \calF^{G}(P_{0})$ be all maximal parabolic subgroups containing $P$.
Let $Q = \bigcap_{i=1}^{k-1}P_{i}$. Let $X \in \Ugo(G)$ and $t > 0$. Let $\{X_{i}\}_{i \in I} \sbs  \Ugo(G)$
be given by Lemma \ref{lem:I210}. We apply this Lemma to the function $\phi_{Q, G}$ 
and the maximal parabolic subgroup $P_{k}$. By induction hypothesis, inequality 
\eqref{eq:I210} holds for all $X_{i}$ and exponent $\la = \mu - t\al_{Q}$. 
Applying the lemma 
and observing that $(\phi_{Q, G})_{P_{k}, G} = \phi_{P,G}$, we get the desired result.
 \edem

We define $\ago_{0'}^{G}$ as the orthogonal complement of $\zgo_{G}$ in $\ago_{0'}^{G}$. 
For $R \in \calF^{G}(P_{0'})$ we define also $\zgo_{R}^{G, +}$ to be $\zgo_{R}^{+} \cap \zgo_{R}^{G}$.
Siegel domains are discussed in Paragraph \ref{ssec:redTh}.
We have then the Siegel domain of $G'(\A)$ with respect to the subgroup $P_{0'}$. Let 
$\dl$ be a subset of $\ago_{0'}^{G}$. We define the following subset of the Siegel set $\Sgl(G')$
\[
\Sgl(G', \dl) = 
\{cak \ | \ c \in \omega', \, k \in K', \, a \in A_{0'}^{\infty}, \, H_{0'}(a) \in  \dl \cap  \ago_{0'}^{+}\}. 
\]
Note that this is a subset of $\Sgl(G')$ since $T_{G'} \in - \ago_{0'}^{+}$.

\btheo\label{thm:intCrit}
Let $\phi \in \calA(G)$ be a form with a central character. 
Then 
\[
\int_{[G']^{G, 1}} |\phi (x) \xi(x)| \,dx < \infty
\]
if and only if
for all $P \in \calF^{G, max}(P_{0}')$ and all $\la \in \calE_{P}(\phi)'$
we have
\[
\Rel(\la + \xi) + \upla_{P} \in - \rint (\zgo_{P}^{G} \cap \brzgop)^{\vee}.
\]
If the above condition holds, we have $\phi \in \calA(G)^{\xi-reg}$ and
\[
\calP(\xi, \phi) = \int_{[G']^{G,1}}\phi(x)\xi(x)\,dx.
\]
\etheo
\bdem
The if part of the theorem is an immediate consequence of the inversion formula \ref{lem:invForm} 
and the Theorem \ref{thm:rapDec}.

Let's prove the necessity. 
We follow closely the proof of the square-integrability criterion of automorphic forms proved 
in \cite{moeWald2}, Lemma I.4.11. 
Fix $P  \in \calF^{G, max}(P_{0}')$. 
The decomposition \eqref{eq:bigFan} in the case $P'= P_{0}'$ shows 
that there exists an $R \in \calF^{G}(P_{0}')$ contained in $P$
such that $\zgo_{R}^{+}$ is open in $\ago_{0'}^{+}$. 
Let's fix it. We will assume that $R,P \in \calF^{G}(P_{0})$, where $P_{0}$ was fixed in the beginning of this paragraph.  
Let $\phi \in \calA(G)$ be such that
\[
\int_{[G']^{G, 1}} |\phi (x) \xi(x)| \,dx < \infty.
\]
We know then that
\begin{equation}\label{eq:intCrit1}
\int_{\Sgl(G', \zgo_{R}^{G, +})} |\phi (cak) \xi(cak)| e^{\langle - 2\rho_{0'}, H_{0'}(a) \rangle}\,dcdadk < \infty.
\end{equation}
From identity \eqref{eq:constAlt} we have
\[
\phi - \phi_{P} = \sum_{P \sbn Q \sbs G}\phi_{P,Q}.
\]
For $Q$ as above, put $\al_{P}^{Q} = \sum_{\al \in \Delta_{0}^{Q} \smin \Delta_{0}^{P}}\al$.
It follows from Corollary \ref{cor:I210} that there exist a  $\mu \in \ago_{0'}$ 
such that 
for all $t > 0$ 
there exists a constant $c$ such that 
for all $x \in \Sgl(G', \zgo_{R}^{G, +})$ we have
\begin{equation}\label{eq:intCrit2}
|\phi - \phi_{P}|(x) \le c e^{\langle \mu, H_{0'}(x) \rangle}\sum_{P \sbn Q \sbs G}
e^{\langle -t\al_{P}^{S}, H_{0'}(x) \rangle}.
\end{equation}
Strictly speaking $\Sgl(G', \zgo_{R}^{G, +})$ is only contained in the Siegel domain of $G$ with respect 
to $P_{0}$  up to multiplication by $K'$ on the right, but this difference only influences 
the constants so clearly we can apply the result here. 

Define for $\varepsilon > 0$ and $P \sbn Q \sbs G$
\[
\dl(P, Q, \varepsilon) = \{H \in \zgo_{R}^{G, +} \  | \  \al_{P}^{Q}(H) \ge \varepsilon \|  p_{P, Q}(H)     \| \, \}
\]
where $p_{P,Q} \in \End(\ago_{0'}^{G})$ is the orthogonal projection 
onto the kernel of $\al_{P}^{Q}$ restricted to $\ago_{0'}^{G}$.
Since $\al_{P}^{Q}$ is strictly positive on $\zgo_{R}^{G, +}$ 
we see that $\dl(P, Q, \varepsilon) \cap \zgo_{R}^{G, +}$ is open in $\zgo_{R}^{G, +}$ for $\varepsilon$ small enough. 

We set then
\[
\dl(P,\varepsilon)  = \bigcap_{P \sbn Q \sbs G}\dl(P, Q, \varepsilon) .
\]
For $\varepsilon$ small enough we have then that 
$\dl(P, \varepsilon) \cap \zgo_{R}^{G, +}$ is open in $\zgo_{R}^{G, +}$. 

From \eqref{eq:intCrit1} and \eqref{eq:intCrit2} we see that 
\begin{equation}\label{eq:intCrit3}
\int_{\Sgl(G', \dl(P,\varepsilon) )} |\phi_{P} (cak) \xi(cak)| e^{\langle - 2\rho_{0'}, H_{0'}(a) \rangle}\,dcdadk < \infty.
\end{equation}
for $\varepsilon$ small enough. 

Note that $\dl(P, Q, \varepsilon)$ are stable by addition and multiplication by positive scalars 
(they form cones, although not polyhedral).
We have that $\al_{P}^{Q}$ are strictly positive on $\zgo_{P}^{G,+}$. It follows that for $\varepsilon$ 
small enough, the set $\dl(P,\varepsilon) $ is stable by translation by the half-line $\zgo_{P}^{G, +}$.
Now we can reason exactly as in part (2) of the proof of Lemma I.4.11 in \cite{moeWald2}. 
For $c \in N'(\A)$ and $a \in Z_{P}^{\infty}$ 
we have $\phi_{P}(ca) = \phi_{P}(ac)$. 
For $g$ fixed, we can write $\phi_{P}(ag)$, where $a \in Z_{P}^{\infty} \cap G'(\A)^{1,G}$ as in \eqref{eq:autDecAdpt}
\[
\sum_{i}q_{i}(H_{0'}(a)_{P})e^{\langle \la_{i} + \rho_{P}, H_{0'}(a)_{P}\rangle}
\]
with $\la_{i} \in \calE_{P}(\phi)'$ and $q_{i} \in \C[\zgo_{P}^{G}]$. The polynomials 
depend on $g$ but the exponents do not. Applying Fubini's theorem 
to the convergent integral \eqref{eq:intCrit3}
we obtain 
\[
\int_{X + \zgo_{P}^{G, +}}
|\sum_{i}q_{i}(H)e^{\langle \Rel(\la_{i} + \xi) + \upla_{P}, H \rangle}|\, dH < \infty
\]
for some $X \in \zgo_{P}^{G, +}$ (that depends on $g$)
for almost all $g$ 
of the form $ca'k$ with $c \in \omega'$, $k \in K'$ and $a' \in A_{0'}^{\infty} \cap M_{P \cap G'}(\A)^{1,P}$
(in the measure theoretic sense). This implies $\Rel(\la_{i} + \xi) + \upla_{P} \in - \rint (\zgo_{P}^{G} \cap \brzgop)^{\vee}$ 
for all $\la_{i} \in \calE_{P}(\phi)'$
as desired.

\edem

\bibliographystyle{abbrv}

\begin{thebibliography}{10}

\bibitem{arthur3}
J.~Arthur.
\newblock A trace formula for reductive groups. {I}. {T}erms associated to
  classes in {$G({\bf Q})$}.
\newblock {\em Duke Mathematical Journal}, 45(4):911--952, 1978.

\bibitem{arthur5}
J.~Arthur.
\newblock A trace formula for reductive groups. {II}. {A}pplications of a
  truncation operator.
\newblock {\em Compositio Math.}, 40(1):87--121, 1980.

\bibitem{arthur2}
J.~Arthur.
\newblock The trace formula in invariant form.
\newblock {\em Ann. of Math. (2)}, 114(1):1--74, 1981.

\bibitem{barvinok}
A.~{Barvinok}.
\newblock {\em {Integer points in polyhedra.}}
\newblock Z\"urich: European Mathematical Society (EMS), 2008.

\bibitem{beuz2}
R.~Beuzart-Plessis.
\newblock Comparison of local spherical characters and the {I}chino-{I}keda
  {C}onjecture for unitary groups.
\newblock {\em Preprint arXiv:1602.06538}, 2016.

\bibitem{brunGub}
W.~Bruns and J.~Gubeladze.
\newblock {\em Polytopes, rings, and {$K$}-theory}.
\newblock Springer Monographs in Mathematics. Springer, Dordrecht, 2009.

\bibitem{cass}
W.~Casselman.
\newblock Extended automorphic forms on the upper half plane.
\newblock {\em Math. Ann.}, 296(4):755--762, 1993.

\bibitem{chaud1}
P.-H. {Chaudouard}.
\newblock {Sur une variante des troncatures d'Arthur}.
\newblock {\em Preprint arXiv:1612.04411}, 2016.

\bibitem{chaudZyd}
P.-H. Chaudouard and M.~Zydor.
\newblock Le transfert singulier pour la formule des traces de jacquet-rallis.
\newblock {\em Preprint arXiv:1611.09656}, 2016.

\bibitem{feiLaOff}
B.~Feigon, E.~Lapid, and O.~Offen.
\newblock On representations distinguished by unitary groups.
\newblock {\em Publ. Math. Inst. Hautes \'{E}tudes Sci.}, 115:185--323, 2012.

\bibitem{gjr1}
D.~Ginzburg, D.~Jiang, and S.~Rallis.
\newblock On the nonvanishing of the central value of the {R}ankin-{S}elberg
  {$L$}-functions.
\newblock {\em J. Amer. Math. Soc.}, 17(3):679--722, 2004.

\bibitem{ginLap}
D.~Ginzburg and E.~Lapid.
\newblock On a conjecture of {J}acquet, {L}ai, and {R}allis: some exceptional
  cases.
\newblock {\em Canad. J. Math.}, 59(6):1323--1340, 2007.

\bibitem{helWan}
A.~G. Helminck and S.~P. Wang.
\newblock On rationality properties of involutions of reductive groups.
\newblock {\em Adv. Math.}, 99(1):26--96, 1993.

\bibitem{ichino1}
A.~Ichino.
\newblock On the regularized {S}iegel-{W}eil formula.
\newblock {\em J. Reine Angew. Math.}, 539:201--234, 2001.

\bibitem{ichYamUn}
A.~Ichino and S.~Yamana.
\newblock Periods of automorphic forms: the case of {$(\rm U_{n+1} \times \rm
  U_{n}, \rm U_{n})$}.
\newblock {\em J. Reine Angew. Math.}
\newblock To appear.

\bibitem{ichYamGl}
A.~Ichino and S.~Yamana.
\newblock Periods of automorphic form: the case of {$(\mathrm{GL}_{n+1} \times
  \mathrm{GL}_{n}, \mathrm{GL}_{n})$}.
\newblock {\em Compositio Math.}, 151(4):665--712, 2015.

\bibitem{jacq1}
H.~Jacquet.
\newblock Sur un r\'{e}sultat de {W}aldspurger.
\newblock {\em Ann. Sci. \'{E}cole Norm. Sup. (4)}, 19(2):185--229, 1986.

\bibitem{jacq2}
H.~Jacquet.
\newblock Sur un r\'esultat de {W}aldspurger. {II}.
\newblock {\em Compositio Math.}, 63(3):315--389, 1987.

\bibitem{jacqLai}
H.~Jacquet and K.~F. Lai.
\newblock A relative trace formula.
\newblock {\em Compositio Math.}, 54(2):243--310, 1985.

\bibitem{jacqLaiRal}
H.~Jacquet, K.~F. Lai, and S.~Rallis.
\newblock A trace formula for symmetric spaces.
\newblock {\em Duke Math. J.}, 70(2):305--372, 1993.

\bibitem{jlr}
H.~Jacquet, E.~Lapid, and J.~Rogawski.
\newblock Periods of automorphic forms.
\newblock {\em J. Amer. Math. Soc.}, 12(1):173--240, 1999.

\bibitem{jacqRall2}
H.~Jacquet and S.~Rallis.
\newblock Symplectic periods.
\newblock {\em J. Reine Angew. Math.}, 423:175--197, 1992.

\bibitem{jacqZag}
H.~Jacquet and D.~Zagier.
\newblock Eisenstein series and the {S}elberg trace formula. {II}.
\newblock {\em Trans. Amer. Math. Soc.}, 300(1):1--48, 1987.

\bibitem{jiang}
D.~Jiang.
\newblock {$G_2$}-periods and residual representations.
\newblock {\em J. Reine Angew. Math.}, 497:17--46, 1998.

\bibitem{kudRal}
S.~S. Kudla and S.~Rallis.
\newblock A regularized {S}iegel-{W}eil formula: the first term identity.
\newblock {\em Ann. of Math. (2)}, 140(1):1--80, 1994.

\bibitem{labWal}
J.-P. Labesse and J.-L. Waldspurger.
\newblock {\em La formule des traces tordue d{'}apr\`es le {F}riday {M}orning
  {S}eminar}.
\newblock Paris, 2009.

\bibitem{lapOff}
E.~Lapid and O.~Offen.
\newblock On the distinguished spectrum of {$Sp(2n)$} with respect to {$Sp(n)
  \times Sp(n)$}.
\newblock {\em Kyoto J. Math.}, 58(1):101--171, 2018.

\bibitem{lapid}
E.~M. Lapid.
\newblock On the fine spectral expansion of {J}acquet's relative trace formula.
\newblock {\em J. Inst. Math. Jussieu}, 5(2):263--308, 2006.

\bibitem{lapRog}
E.~M. Lapid and J.~D. Rogawski.
\newblock Periods of {E}isenstein series: the {G}alois case.
\newblock {\em Duke Math. J.}, 120(1):153--226, 2003.

\bibitem{lawrence}
J.~Lawrence.
\newblock Rational-function-valued valuations on polyhedra.
\newblock In {\em Discrete and computational geometry ({N}ew {B}runswick, {NJ},
  1989/1990)}, volume~6 of {\em DIMACS Ser. Discrete Math. Theoret. Comput.
  Sci.}, pages 199--208. Amer. Math. Soc., Providence, RI, 1991.

\bibitem{levy}
J.~Levy.
\newblock A truncated integral of the {P}oisson summation formula.
\newblock {\em Canadian Journal of Mathematics}, 53(1):122--160, 2001.

\bibitem{miVen}
P.~Michel and A.~Venkatesh.
\newblock The subconvexity problem for {${\rm GL}_2$}.
\newblock {\em Publ. Math. Inst. Hautes \'{E}tudes Sci.}, (111):171--271, 2010.

\bibitem{moeWald2}
C.~M{\oe}glin and J.-L. Waldspurger.
\newblock {\em Spectral decomposition and {E}isenstein series}, volume 113 of
  {\em Cambridge Tracts in Mathematics}.
\newblock Cambridge University Press, Cambridge, 1995.
\newblock Une paraphrase de l'{\'E}criture [A paraphrase of Scripture].

\bibitem{offen2}
O.~Offen.
\newblock On symplectic periods of the discrete spectrum of {${\rm GL}_{2n}$}.
\newblock {\em Israel J. Math.}, 154:253--298, 2006.

\bibitem{offen1}
O.~Offen.
\newblock Residual spectrum of {${\rm GL}_{2n}$} distinguished by the
  symplectic group.
\newblock {\em Duke Math. J.}, 134(2):313--357, 2006.

\bibitem{polWanZyd}
A.~Pollack, C.~Wan, and M.~Zydor.
\newblock {A $\mathrm{G}_2$-period of a Fourier coefficient of an Eisenstein
  series on $\mathrm{E}_6$}.
\newblock {\em Israel J. Math.}
\newblock Accepted.

\bibitem{polWanZyd2}
A.~Pollack, C.~Wan, and M.~Zydor.
\newblock {On the residue method for period integrals}.
\newblock {\em Preprint arXiv:1903.02544}, 2019.

\bibitem{khovPukh}
A.~V. Pukhlikov and A.~G. Khovanski\u{\i}.
\newblock The {R}iemann-{R}och theorem for integrals and sums of
  quasipolynomials on virtual polytopes.
\newblock {\em Algebra i Analiz}, 4(4):188--216, 1992.

\bibitem{sak1}
Y.~Sakellaridis.
\newblock The {S}chwartz space of a smooth semi-algebraic stack.
\newblock {\em Selecta Math. (N.S.)}, 22(4):2401--2490, 2016.

\bibitem{sakVen}
Y.~Sakellaridis and A.~Venkatesh.
\newblock Periods and harmonic analysis on spherical varieties.
\newblock {\em Ast\'{e}risque}, (396):viii+360, 2017.

\bibitem{schiff}
O.~Schiffmann.
\newblock Indecomposable vector bundles and stable {H}iggs bundles over smooth
  projective curves.
\newblock {\em Ann. of Math. (2)}, 183(1):297--362, 2016.

\bibitem{schneider}
R.~Schneider.
\newblock Combinatorial identities for polyhedral cones.
\newblock {\em Algebra i Analiz}, 29(1):279--295, 2017.

\bibitem{springer}
T.~A. Springer.
\newblock {\em Linear algebraic groups}, volume~9 of {\em Progress in
  Mathematics}.
\newblock Birkh\"{a}user Boston, Inc., Boston, MA, second edition, 1998.

\bibitem{wu}
H.~{Wu}.
\newblock Deducing selberg trace formula via rankin-selberg method for {${\rm
  GL}_2$}.
\newblock {\em Preprint arXiv:1810.09437}, 2018.

\bibitem{yamana1}
S.~Yamana.
\newblock Symplectic periods of the continuous spectrum of {${\rm GL}(2n)$}.
\newblock {\em Ann. Inst. Fourier (Grenoble)}, 64(4):1561--1580, 2014.

\bibitem{yamana2}
S.~Yamana.
\newblock Periods of residual automorphic forms.
\newblock {\em J. Funct. Anal.}, 268(5):1078--1104, 2015.

\bibitem{zagier}
D.~Zagier.
\newblock The {R}ankin-{S}elberg method for automorphic functions which are not
  of rapid decay.
\newblock {\em J. Fac. Sci. Univ. Tokyo Sect. IA Math.}, 28(3):415--437 (1982),
  1981.

\bibitem{zhengZyd}
H.~Zheng and M.~Zydor.
\newblock A new valuation on polyhedral cones.
\newblock {\em Preprint arXiv:1901.07668}, 2019.

\bibitem{leMoi3}
M.~Zydor.
\newblock Les formules des traces relatives de {J}acquet-{R}allis grossi\`eres.
\newblock {\em J. Reine Angew. Math.}
\newblock Accepted.

\bibitem{leMoi}
M.~Zydor.
\newblock La variante infinit\'{e}simale de la formule des traces de
  {J}acquet-{R}allis pour les groupes unitaires.
\newblock {\em Canad. J. Math.}, 68(6):1382--1435, 2016.

\bibitem{leMoi2}
M.~Zydor.
\newblock La variante infinit\'{e}simale de la formule des traces de
  {J}acquet-{R}allis pour les groupes lin\'{e}aires.
\newblock {\em J. Inst. Math. Jussieu}, 17(4):735--783, 2018.

\end{thebibliography}

\end{document}